\documentclass[10pt,leqno]{article}
\usepackage{amsmath, amssymb, amsfonts, amsthm}
\usepackage{xspace, amscd, graphics}

\newcommand{\Q}{\mathbb Q}
\newcommand{\Z}{\mathbb Z}
\newcommand{\R}{\mathbb R}
\newcommand{\C}{\mathbb C}
\newcommand{\G}{\mathfrak G}

\newcommand{\h}{\mathcal H}
\newcommand{\A}{\mathbb A}

\newcommand{\Oo}{\mathfrak O}


\newcounter{mysubsection}[section]
\title{Algebraic cycles on Hilbert modular fourfolds and poles of $L$-functions}

\author{Dinakar Ramakrishnan\\253-37 Caltech,
Pasadena, CA 91125}
\date{}
\begin{document}
\maketitle

\begin{flushright}To M.S.~Raghunathan on his sixtieth birthday, with admiration
\end{flushright}
\bigskip

\pagestyle{myheadings}

\tableofcontents

\bigskip

\addcontentsline{toc}{section}{Introduction}
\section*{\bf Introduction}

\medskip

Let $K$ be a quartic, totally real number field with ring of
integers $\Oo_K$ and canonical embedding $K \hookrightarrow \R^4$,
given by $x \to (x_\sigma)_{\sigma \in {\rm Hom}(K,\R)}$. Any
congruence subgroup $\Gamma$ of SL$(2, \Oo_K)$ is a discrete
subgroup of ${\rm SL}(2,\R) \times {\rm SL}(2,\R) \times {\rm
SL}(2,\R) \times {\rm SL}(2,\R)$ with finite covolume. It acts
properly discontinuously on the 4-fold product of the upper half
plane $\mathcal H$, with the quotient $\Gamma\backslash{\mathcal
H}^4$ being a complex analytic variety of dimension $4$ with at
most quotient singularities. Let $\A = \R \times \A_f$ denote the
adele ring of $\Q$, $C$ a compact open subgroup of $G(\A_f)$, and
$Y = Y_C$ the corresponding Hilbert modular fourfold, which is a
quasi-projective, normal complex variety, smooth for $C$ small
enough. It comes equipped with the Baily-Borel compactification
$Y^\ast$ and a smooth toroidal compactification $X: = \tilde Y$,
all defined over $\Q$. It is well known that $Y$ is the moduli
space of abelian fourfolds $A$ with level $C$-structure and with
End$(A) \otimes \Q \supset K$. The components of $Y(\C)$ are all
of the form $\Gamma\backslash{\mathcal H}^4$, defined over an
abelian extension of $\Q$. We will consider the cohomology of $X$
in various {\it avatars} like the singular (Betti) version
$H^\ast_B(X(\C), \Q)$ and the \'etale version
$H^\ast_{et}(X_{\overline \Q}, \Q_\ell)$.

For every $q \geq 0$ and for any field $k \supset \Q$, denote by
${\mathcal Z}^q(X_k)$ the $\Q$-vector space generated by the
$k$-rational algebraic cycles of codimension $q$ on $X_{\overline
\Q}$ modulo homological equivalence, by which we mean the
equivalence in $H^{2q}_B(X^\ast(\C), \Q(j))$. If $\overline k$
denotes an algebraic closure of $k$, the absolute Galois group
$\G_k = {\rm Gal}(\overline k/k)$ acts on ${\mathcal
Z}^q(X_{\overline k})$ with ${\mathcal Z}^q(X_k)$ being the
$\Q$-subspace of fixed points. It is well known that ${\mathcal
Z}^0(X_{\overline k}) \simeq \Q$, ${\mathcal Z}^q(X_{\overline k})
\simeq {\mathcal Z}^{4-q}(X_{\overline k})$ and ${\mathcal
Z}^1(X_k)$ is NS$(X_k) \otimes \Q$, where NS$(X_k)$ denotes the
N\'eron-Severi group of $k$-rational divisors on $X$ modulo
algebraic equivalence. The really mysterious one is the group
${\mathcal Z}^2(X_k)$ and there is a dearth of concrete results
about codimension $2$ cycles which are not generated by
intersections of divisors.

Fix a prime $\ell$ and for every $j \geq 0$ denote by
$V_\ell^{(j)}$ the $\ell$-adic Gal$(\overline \Q/\Q)$-
representation defined by $H^j_{et}(X_{\overline \Q}, \Q_\ell)$.
For any {\it number field} $k$, the restriction of $V_\ell^{(j)}$
to $\G_k$ defines an $L$-function
$$
L(s, V_\ell^{(j)}/k) \, = \, \prod\limits_{P} \, \frac{1}{{\rm
det}(I - (NP)^{-s} Fr_P\vert {V_\ell^{(j)}}^{I_P})}
$$
where $P$ runs over the primes in $\Oo_k$ with norm $NP$,
(geometric) Frobenius $Fr_P$ and inertia group $I_P$. There exists
a finite set $S_{k,\ell}$ of primes $P$ containing the primes
above $\ell$ outside which $I_P$ acts trivially on $V_\ell^{(j)}$,
and moreover, by Deligne's proof of the Weil conjectures, the
eigenvalues of $Fr_P$ on $V_\ell^{(j)}$ for any $P$ outside
$S_{k,\ell}$ are of absolute value $(NP)^{j/2}$. This shows that
for any finite set $S$ of primes containing $S_{k,\ell}$, the {\it
incomplete $L$-function} $L^S(s,V_\ell^{(j)}/k)$, defined as the
Euler product above except with the factors at $S$ removed, is
absolutely convergent in $\Re(s) > \frac{j}{2}+1$. It is not {\it
a priori} clear that this $L$-function is even defined at the {\it
Tate point} $s = \frac{j}{2}+1$. When $j=2m$, a celebrated {\it
conjecture of Tate} says that the $L$-function is meromorphic at
this point and moreover, that the {\it order of pole} there is
equal to the dimension of ${\mathcal Z}^m(X_k)$. The object of
this paper is to provide some non-trivial evidence for it for $m =
2$. Since we will be interested only in the cohomology in degree
$4$, we will write $V_\ell$ for $V_\ell^{(4)}$. For any non-zero
ideal $\mathfrak N$ in $\Oo_K$, let $C_0({\mathfrak N})$, resp.
$C_1({\mathfrak N})$, denote the compact open subgroup of $G_f$
given as the product over the finite places $v$ of
$C_{0,v}(P_v^{v({\mathfrak N})})$, resp.
$C_{1,v}(P_v^{v({\mathfrak N})})$; here $C_{0,v}(P_v^{r})$ denotes
the group of matrices $\begin{pmatrix} a & b\\ c & d\end{pmatrix}$
in GL$(2, \Oo_{K_v})$ with $c \equiv 0 \, (\bmod P_v^r)$, and
$C_{1,v}(P_v^{r})$ signifies the subgroup of such matrices with $d
\equiv 1 \, (\bmod P_v^r)$ . (By convention, $P_v^0$ is
$\Oo_{K_v}$ and so there is no condition when $r=0$.) Now let
$\Psi$ be any weight $1$ Hecke character of a totally imaginary
quadratic extension $M$ of $K$. Then it is well known that $\Psi$
defines a Hilbert modular newform $g_\Psi$, contributing to
$V_\ell$ relative to a level subgroup $C$ containing
$C_1({\mathfrak N})$ if ${\mathfrak N}$ is divisible by the
relative discriminant of $M/K$ times the norm of the conductor of
$\Psi$. By the {\it non-CM part} of $V_\ell$ we will mean the
quotient $V_\ell'$ of $V_\ell$ by the space spanned by the classes
attached to such Hecke characters. Let ${\mathcal Z}^2(X_k)'$
denote, for any number field $k$, the dimension of the non-CM part
of ${\mathcal Z}^2(X_k)$. Denote also by $Ta_{\ell, k}(X)$ the
space of Tate classes in $V_\ell(2)$ (see section 1 for a
definition), and by $Ta_{\ell, k}(X)'$ its non-CM part.

Throughout the paper we will use the classical term {\it abelian
number field} to mean a finite abelian  extension of $\Q$.

\medskip

\noindent{\bf Theorem A} \, \it Let $K$ be a quartic, totally real
number field containing a quadratic subfield, $C$ a compact open
subgroup of $G(\A_f)$, $X$ a Hilbert modular fourfold of level
$C$, and $S$ a finite set of primes containing $S_{k,\ell}$. Then
\begin{enumerate}
\item[(i)] \, $L^S(s, V_\ell/k)$ extends, for any abelian number field $k$,
to a meromorphic function
on all of $\C$ and satisfies, after being multiplied by suitable
local factors at $S$ and at the archimedean places, a functional
equation relating $s$ to $5-s$;
\item[(ii)] \, Let $K/\Q$ be
Galois. Then we have
$$
{\rm dim}_\Q {\mathcal Z}^2(X_k)' \, = \, -{\rm ord}_{s=3} L^S(s,
V_\ell'/k),
$$
for any abelian number field $k$, and
$$
{\rm dim}_\Q {\mathcal Z}^2(X_k)' \, = \, {\rm dim}_{\Q_\ell}
Ta_{\ell, k}(X)'.
$$
for any number field $k$.
\end{enumerate}
\rm

\medskip

The proof will show that there are algebraic cycles of codimension
$2$ which are not in the $\Q$-linear span of intersections of
divisor classes. When $K/\Q$ is non-Galois, it can be shown that
the dimension of ${\mathcal Z}^2(X_k)$ is still bounded, for $k$
abelian, by the order of pole of $L^S(s, V_\ell)$ at $s=3$, and
the desired equality will follow from the proof in some cases.

\medskip

Here are some general philosophical remarks. Given a general
smooth, projective variety $X$, it is hard to find {\it any}
collection of explicit algebraic cycles on it. For Shimura
varieties one is fortunate in this regard as there are some
natural cycles $Z$, though not nearly enough, defined by {\it
Shimura subvarieties} and their translates by Hecke
correspondences, as well as certain twists. There are also cycles
supported on certain rational varieties, and their intersections,
occurring in the smooth resolution at infinity. One can often
write down explicitly a basis for the $(p,p)$-cohomology for any
$p$. Still, given a differential form $\omega$ representing such a
cohomology class and a dimension $p$ algebraic cycle $Z$, how can
one show that the integral of $\omega$ over $Z$ is non-zero? One
way, and this is the tack we take in this paper, is to express
$\int_Z \omega$ as the residue at a pole of some $L$-function
which is non-zero for some reason, for example because of the
knowledge of the exact order of pole at that point. The {\it
novelty}, if any, in the situation considered here, is that the
$L$-functions which appear this way are {\it not} the
$L$-functions of the varieties involved, but rather certain
associated ones, and luckily there is a (partial) coincidence of
poles for the different $L$-functions. In certain situations we
can also check the Hodge conjecture as seen in the result below:

\medskip

\noindent{\bf Theorem B} \, \it Let $K$ be a quartic, totally
real, Galois number field, and let $X$ (resp. $X(1)$) be a smooth,
projective, Hilbert modular fourfold associated to the level
subgroup $C_0({\mathfrak N})$ for some square-free ideal
${\mathfrak N} \ne \mathfrak O_K$ (resp. for ${\mathfrak N} =
\mathfrak O_K$). Then the Hodge cycles on $X_\C$ in codimension
$2$, which are not pull-backs of classes in $X(1)_\C$, are all
algebraic. For various $\mathfrak N$, there are such algebraic
cycles $Z$ of codimension $2$ on $X$ which are not intersections
of divisors. \rm

\medskip

For a statement of the Hodge conjecture (and the definition of
Hodge cycles), see section 1. The proof of Theorem B will be given
in section 10. However, it should be acknowledged that the proof
given there depends partially on a joint result with V.~Kumar
Murty on a comparison of rational Hodge structures (see Theorem
10.4) which we plan to publish elsewhere in a more general
setting. If this is not satisfactory to some readers, they can
focus just on the proof of Theorem A which is totally
self-contained.

\medskip

This paper provides an extension of the work of Harder, Langlands
and Rapoport ([HLR]) on the Tate conjecture for {\it divisors} on
Hilbert modular {\it surfaces} over abelian number fields. This
generalization is {\it not routine}, however, due to certain
subtle problems, the least of which is that we are working with
codimension $2$ cycles. It is perhaps helpful to elaborate.

The {\it first difficulty}, which is analytic, is to show that
$L(s, V_\ell)$ is meromorphic at $s=3$, which we will discuss
below. Granting that, the {\it next step} is to look for
$\Q$-rational algebraic cycles $Z$ on $\tilde X$ not entirely
supported on $X^\infty$, with a view to matching their
(homological) non-triviality with the possible poles of $L(s,
V_\ell)$. Under the hypothesis that $K$ contains a quadratic
subfield $F$, the natural cycles to consider are Hecke translates
of the Hilbert modular surface defined by GL$(2)/F$; these are the
analogs of the Hirzebruch-Zagier cycles investigated in [HLR], and
earlier -- in a concrete form -- in [HZ]. It turns out that their
non-triviality is linked to the poles of the {\it Asai}
$L$-functions $L(s, As_{K/F})$ associated to $K/F$, which are of
degree $4$ over $F$ and $8$ over $\Q$. More precisely, if $\pi$ is
a cusp form on GL$(2)/K$ of weight $2$, then the $\pi$-part of
$L(s, As_{K/F})$, written $L(s, \pi; r_{K/F})$, has a pole at the
right edge iff the {\bf period integral}
$$
\int\limits_{{\rm GL}(2,F)Z(\A_F)\backslash {\rm GL}(2, \A_F)} \,
\phi(g) dg
$$
is non-zero for some function $\phi$ in the space of $\pi$, where
$Z$ denotes the center of GL$(2)$. This period is simply the
residue at the edge of convergence ({\it Tate point}) of $L(s,
\pi; r_{K/F})$, which has a representation as the integral over
${\rm GL}(2,F)Z(\A_F)\backslash {\rm GL}(2, \A_F)$ of $\phi$ times
an Eisenstein series $E(s)$ on GL$(2)/F$.

The {\it second difficulty} is that these functions $L(s, \pi;
r_{K/F})$ do {\bf not} divide $L(s, V_\ell)$. Fortunately for us
we are able to show, and this is a key point, that the poles of
$L(s, As_{K/F})$ give rise, over abelian fields $k$, to poles of
$L(s, V_\ell)$. But they do not account for all the poles of $L(s,
V_\ell)$. When $K/\Q$ is Galois and $\pi$ non-CM, however, we are
able to show that all such poles are accounted for by suitable
twists of (Hecke translates of) embeddings of Hilbert modular
surfaces relative to {\it all} the quadratic fields contained in
$K$.

The {\it third difficulty} comes up in the biquadratic case when
the order of pole of $L(s, \pi; r_{K/F})$ can be $2$ for certain
cusp forms $\pi$. Then it is not enough to consider only the
(Hecke translates of) cycles coming from one quadratic
subextension, and more importantly, when we consider a pair of
cycles coming from two different quadratic subextensions, a key
point of the proof is to show that these two, together with a
suitable twist of one of them (see section 9), must span a plane
in the $\pi_f$-component of the homology. The referee has
indicated an alternate, elegant way to handle this, which has been
described in Remark 9.17, but we have left our proof intact as we
think this method could be of use in other situations.

\medskip

One knows (cf. [BrL], [La1]) that the main part of $L(s, V_\ell)$
is a product, over cusp forms $\pi$ of GL$(2)/K$ of weight $2$, of
certain Asai $L$-functions $L(s,\pi; r_{K/\Q})$ of degree $16$.
Under the hypothesis that $K$ contains a quadratic subfield $F$ we
establish (see section 7 below) the meromorphic continuation and
functional equation for such $L$-functions. We cannot analyze
their poles except, luckily, at the right edge of absolute
convergence, which is what is relevant for the Tate conjecture. To
be precise, for any such $\pi$, and for any Dirichlet character
$\nu$ of $\Q$, the $\nu$-twisted $L$-function $L(s, \pi; r_{K/\Q}
\otimes \nu)$ turns out to admit a pole at the right edge iff a
suitable twist of $\pi$ is a base change from a quadratic subfield
of $K$; this is as in the case of Hilbert modular surfaces
([HLR]). But more interestingly, we show that such an $L$-function
has a {\it double pole} at the Tate point iff a twist of $\pi$ is
a base change all the way from $\Q$ {\it and} $K/\Q$ is
biquadratic.

\medskip

The CM case is not treated in this paper, partly to not make this
paper longer, and partly because there are more Tate classes than
we can handle in certain situations. In other words, there are
exotic ones which cannot be accounted for by the analogs of
Hirzebruch-Zagier cycles (and their twists), and such classes can
exist even over abelian fields -- see Remark 9.18; this was {\it
not} the case for Hilbert modular surfaces. It is likely that a
variant of the method of [Mu-Ra] can match such Tate cycles with
corresponding Hodge cycles, but since the Hodge conjecture is not
known in codimension $2$, the analysis grinds to a halt.

\medskip

Quite generally, given any finite extension $K/F$ of number fields
of degree $d$, and any isobaric automorphic form $\pi$ on
GL$(n)/K$, one can associate an Asai $L$-function, denoted $L(s,
\pi; r_{K/F})$ (see section 6 for a definition), of degree $n^d$.
This is an analogue of tensor induction on the Galois side. For
example, suppose $K/F$ is cyclic and assume the existence of an
$n$-dimensional representation $\sigma$ (over $\C$ or $\overline
\Q_\ell$) of Gal$(\overline K/K)$ with the same $L$-function as
$\pi$. For every $g \in G: = {\rm Hom}(K, \overline K)$, let
$\sigma^{[g]}$ denote the representation of Gal$(\overline K/K)$
defined by $\beta \to \sigma(g\beta g^{-1})$. Then the
representation $\otimes_{g \in G} \, \sigma^{[g]}$ is
$G$-invariant and extends to a non-unique $n^d$-dimensional
representation of Gal$(\overline F/F)$. One can define a
particular choice of an extension, called the {\it Asai
representation} and denoted $As(\sigma)$. The $L$-function of
$As(\sigma)$ should be the same as $L(s, \pi; r_{K/F})$. The {\it
principle of functoriality} implies that the Asai $L$-function,
whether or not $\pi$ is associated to any Galois representation,
should have meromorphic continuation and a standard functional
equation, and this is known in the $(n,2)$-case for arbitrary $n$
(see section 6 below). Classically, Asai established the requisite
properties, and in fact the location of all the poles, for the
case $(n,d) = (2,2)$ when $F=\Q$, $K$ real quadratic and $\pi$
holomorphic (Hilbert modular) newform over $K$; the general
$(2,2)$-case was treated in [HLR], and then in [JY] by using the
{\it relative trace formula}. The case $(n,d) = (2,3)$ was treated
in ([RPS]), providing a non-trivial extension of Garrett's
construction for the triple product $L$-function, and a complete
location of the poles was then given in [Ik]. The meromorphic
continuation and functional equation follows from the work of
Shahidi ([Sh])in the $(n,2)$-case for any $n$, and the possible
poles at the right edge were analyzed in [F$\ell$]. The main thing
for us is that we can also treat, using [Ra3], the $(2,4)$-case
under the hypothesis that $K$ contains a quadratic subextension
$F$.

For any $(n,d)$, the principle of functoriality implies moreover
that there is an isobaric automorphic form $As_{K/F}(\pi)$ on
GL$(n^d)/F$ whose standard $L$-function coincides with $L(s, \pi;
r_{K/F})$. When $(n,d)=(2,2)$ this was established in [Ra2], and
we know of no such result in this direction when $n > 2$ and $d >
1$. This modularity result is crucially used below in deducing the
needed properties of $L(s, \pi; r_{K/\Q})$.

\medskip

I lectured on this material during July 2002 at the conference on
Automorphic Forms at Park City, UT, and then later at Columbia
University, NY. I benefited from the feedback and interest I
received from various people at those places including L.~Clozel,
H.~Jacquet, L.~Saper, and S.~Zhang. I would also like to thank
Jacob Murre for a helpful conversation, Arvind Nair for pointing
out the need for elaboration at a point in section 3, and finally
the referee for a careful reading of the paper and for making
several good comments.

\medskip

The last section contains some remarks on my earlier papers. The
main thing here is to give some details, sought by Joe Shalika, of
the proof of the key Lemma 3.4.9 of [Ra2] giving Sobolev
inequalities for eigenfunctions of the Casimir operator. It also
contains some clarifications/refinements sought by E.~Lapid and
M.~Krishnamurthy on [Ra2] and [Ra3] respectively, as well as some
errata for [Ra2,3].

\medskip

When I began my graduate studies in Mathematics at Columbia
University in the Fall of 1975, I had the good fortune to attend
Raghunathan's exciting, challenging and (very) fast course on
Arithmetic groups. I learnt a lot from that one semester course,
which has come in handy at various times in my own work which
tends to traverse nearby fields like Automorphic Forms and
Arithmetical Geometry. Raghunathan has also been a very
encouraging and friendly figure over the years. It is a great
pleasure to dedicate this article to him. Thanks are also due to
the NSF for financial support through the grant DMS--0100372.

\vskip 0.2in

\section{The algebraic versus analytic rank}

\bigskip

Denote by $\overline \Q$ be the algebraic closure of $\Q$ in $\C$.
Let $X$ be any smooth projective variety of dimension $d$ over
$\Q$. A $\overline \Q$-rational {\bf algebraic cycle of
codimension $q \leq d$} is a finite formal sum
$$
Z = \sum_{i=1}^m \alpha_i Z_i,
\leqno(1.1)
$$
with each $\alpha_i \in \Q$ and $Z_i$ a closed, irreducible
subvariety of $X_{\overline \Q} = X \times_\Q \overline \Q$ of
codimension $q$. Evidently the collection of all such $Z$ forms a
$\Q$-vector space, denoted $C^q(X_{\overline \Q})$.

Integration of differential forms of degree $2d-2q$ over such a
cycle $Z$ defines a class
$$
[Z] \, \in \, H_{2d-2q}(X(\C), \Q) \, \simeq \, H^{2q}_B(X(\C),
\Q)(q),
\leqno(1.2)
$$
where the subscript $B$ signifies taking the {\it singular} (or
Betti) cohomology, and $(q)$ denotes twisting by $\Q(q) = (2\pi
i)^q\Q$. Since $\Q(q)$ is the unique Hodge structure of rank $1$,
weight $-2q$ and bidegree $(-q,-q)$, $H^{2q}_B(X(\C), \Q)(q)$ is a
Hodge structure of weight $0$. Recall that there is a Hodge
decomposition
$$
H^{2q}(X(\C), \C) \, = \, \oplus_{i+j=2q} H^{i,j}(X),
\leqno(1.3)
$$
where the classes in $H^{i,j}(X)$ are represented by differential
forms of bidegree $(i,j)$.

A basic fact is that $[Z]$ lies in the subspace of {\bf Hodge
cycles} of codimension $q$ on $X(\C)$, defined to be
$$
Hg^q(X) \, = \, (H^{2q}_B(X(\C), \Q) \cap H^{q,q}(X))(q).
\leqno(1.4)
$$
The {\bf Hodge conjecture} says that every Hodge cycle on $X(\C)$
is of the form $[Z]$ for an algebraic cycle $Z$; it is known to
hold, by Lefschetz, for $q=1$.

\medskip

Put
$$
Z \, \equiv \, 0 \, \Leftrightarrow \, [Z] \, = \, 0.
\leqno(1.5)
$$
This gives an equivalence, called the {\bf homological
equivalence}, for algebraic cycles. Set
$$
{\mathcal Z}^q(X_{\overline \Q}) \, = \, C^q(X_{\overline \Q})/
\equiv.
\leqno(1.6)
$$

\medskip

Since $X$ is defined over $\Q$, the absolute Galois group
$\mathcal G_\Q = $Gal$(\overline \Q/\Q)$ permutes the closed
irreducible subvarieties of $X_{\overline \Q}$. Hence we get a
Galois action $(\sigma, Z) \to Z^\sigma$ on $C^q(X_{\overline
\Q})$ given by
$$
Z^\sigma = \sum_{i=1}^m \alpha_i Z_i^\sigma \quad if \quad  Z =
\sum_{i=1}^m \alpha_i Z_i. \leqno(1.7)
$$

Now fix a prime $\ell$ and consider the $\ell$-adic cohomology
groups of $X_{\overline \Q}$, on which there is a natural action
of $\mathcal G_\Q$. Then we also have $\ell$-adic cycle classes
$[Z]_\ell$ in $H^{2q}_{et}(X_{\overline \Q}, \Q_\ell)(q)$, with
$(q)$ denoting the tensoring with the Galois representation
$\Q_\ell(q) = \Z_\ell(q) \otimes_{\Z_\ell} \Q_\ell$, where
$Z_\ell(q)$ denotes the inverse limit $\lim\limits_{n}
\mu_{\ell^n}^{\otimes q}$ and $\mu_{\ell^n}$ the group of
$\ell^n$-th roots of unity in $\overline \Q$. One knows the
following:
\begin{enumerate}
\item[] There is an isomorphism of $\Q_\ell$-vector spaces
$$
H^{2q}_B(X(\C), \Q)(q) \otimes_\Q \Q_\ell \, \simeq \,
H^{2q}_{et}(X_{\overline \Q}, \Q_\ell)(q) \leqno(1.8)
$$
such that the image of $[Z]$ is $[Z]_\ell$;
\item[] The
$\ell$-adic cycle class map $Z \to [Z]_\ell$ is Galois
equivariant, i.e.,
$$
[Z^\sigma]_\ell \, = \, [Z]_\ell^\sigma,  \, \, \forall \, \sigma
\in {\mathcal G}_\Q. \leqno(1.9)
$$
\end{enumerate}
It follows that if $Z \equiv 0$, then $Z^\sigma \equiv 0$ for all
$\sigma$, and this gives us an action of $\mathcal G_\Q$ on
${\mathcal Z}^q(X_{\overline \Q})$. For any field $k \subset
\overline \Q$, we define the {\bf group of $k$-rational algebraic
cycles} of codimension $q$ {\it modulo homological equivalence} to
be
$$
{\mathcal Z}^q_k(X) : \, = \, {\mathcal Z}^q(X_{\overline
\Q})^{\mathcal G_k},
\leqno(1.10)
$$
which is finite dimensional because it identifies with a
$\Q$-subspace of $H^{2q}(X(\C), \Q)(q)$.

The {\bf algebraic rank} of $X$ over $k$ {\it in codimension $q$}
is defined to be
$$
r_{{\rm alg}, k}^q(X) \, = \, {\rm dim}_\Q {\mathcal Z}^q_k(X).
\leqno(1.11)
$$

\medskip

For any $j \leq 2d$, put
$$
V_\ell^{(j)} \, = \, H^j_{et}(X_{\overline \Q}, \Q_\ell),
\leqno(1.12)
$$
which is a finite dimensional, continuous representation of
$\mathcal G_\Q$. Let $S_k$ be a finite set of primes $P$ in $k$
such that either $P \mid\ell$ or $V_\ell^{(2q)}$ is ramified at
$P$. The incomplete {\bf $L$-function} attached to $(V_\ell^{(j)},
k)$ is
$$
L^S(s, V_\ell^{(j)}/k) \, = \, \prod\limits_{P \notin S_k} {\rm
det}(I-Fr_P P^{-s} \vert V_\ell^{(j)})^{-1},
\leqno(1.13)
$$
which converges absolutely in the right half plane $\{\Re(s) >
j/2+1\}$ by Deligne's proof of the Weil conjectures, which assert
that the inverse roots of the geometric Frobenius elements $Fr_P$
are all of absolute value $N(P)^{j/2}$.

Take $j = 2q$. One expects $L^S(s, V_\ell^{(2q)})$ to admit a
meromorphic continuation and satisfy, with suitable factors at
$S_k$ and infinity, a functional equation relating $s$ and
$2q+1-s$. All we need for the Tate conjecture, however, is that
this $L$-function is meromorphic at the {\it Tate point}, namely
the {\it edge of convergence} $s=q+1$. Admitting this leads to the
following definition, for every number field $k$, of the {\bf
analytic rank} of $X$ over $k$ {\it in codimension $q$}:
$$
r_{{\rm an}, k}^q(X) \, = \, -{\rm ord}_{s=q+1} L^S(s,
V_{\ell}^{(2q)}/k),
\leqno(1.14)
$$
where $V_{\ell, k}$ denotes the restriction of $V_\ell$ to
$\mathcal G_k$. We have

\medskip

\noindent{\bf Conjecture I} (Tate) \, \it
$$
r_{{\rm alg}, k}^q(X) \, = \, r_{{\rm an},k}^q(X).
\leqno(1.15)
$$
\rm

\medskip

This is true for $q = 0, d$, and there is some positive evidence
for divisors ($q=1$).

\medskip

Finally, by the Galois equivariance of the $\ell$-adic cycle class
map, $[Z]_\ell$ lies in the space of {\bf Tate cycles} over $k$:
$$
Ta_{\ell, k}^q(X) : = \, H^{2q}_{et}(X_{\overline \Q},
\Q_\ell(q))^{\mathcal G_k}, \leqno(1.16)
$$
for all $Z \in {\mathcal Z}^q_k(X)$. The {\bf $\ell$-adic cycle
rank} of $X$ over $k$ in {\it codimension $q$} is then give by
$$
r_{\ell,k}^q(X) \, = \, {\rm dim}_{\Q_\ell} Ta_{\ell, k}^q(X).
\leqno(1.17)
$$
One also has the following

\medskip

\noindent{\bf Conjecture II} (Tate) \, \it
$$
r_{\ell, k}^q(X) \, = \, r_{{\rm an}, k}^q(X).
\leqno(1.18)
$$
\rm

\medskip

Since ${\mathcal Z}_k^q(X)$ injects into $Ta_{\ell, k}^q(X)$, we
get the following

\medskip

\noindent{\bf Proposition 1.19} \, \it
$$
r_{{\rm alg}, k}^q(X) \, \leq \, r_{\ell, k}^q(X).
$$
\rm

\medskip

The existence of the {\bf Hodge-Tate decomposition} ([F]) implies
the following useful bound (where $h^{q,q}(X)$ denotes the
dimension of $H^{q,q}(X_\C)$):

\medskip

\noindent{\bf Proposition 1.20} \, \it
$$
r_{\ell, k}^q(X) \, \leq \, h^{q,q}(X).
$$
\rm

\medskip

Let $r_{Hg}^q(X)$ denotes the dimension of the space $Hg^q(X_\C)$
of {\bf Hodge cycles} (of codimension $q$) on $X_\C$. We trivially
have
$$
r_{Hg}^q(X) \, \leq \, h^{q,q}(X), \leqno(1.21)
$$
and the Hodge conjecture is the statement
$$
r_{{\rm alg}}^q(X) \, = \, r_{Hg}^q(X_\C),
$$
where $r_{\rm alg}^q(X)$ denotes the dimension of algebraic cycles
in $Hg^q(X_\C)$, which, by the comparison theorem (relating the
Betti and \'etale cohomology) and the proper base change theorem
(relating the \'etale cohomology of $X_{\overline \Q}$ with that
of $X_\C$, identifies with $r_{{\rm alg}, k}(X)$ for $k$ a large
enough number field.

\medskip

If $X$ is an abelian variety over a number field $k$, it is known
by Deligne ([DMOS]) that for any prime $\ell$,
$$
r_{Hg}^q(X) \, \leq \, r_{\ell}^q(X_{\overline k})
\leqno(1.22)
$$
for $k$ large enough. This will be used in section 10.

\vskip 0.2in

\section{Hilbert modular varieties}

\bigskip

Let $K$ be a totally real number field of degree $d$, with ring of
integers $\Oo_K$ and adele ring $\A_K = K_\infty \times \A_{K,f}$,
where the ring $\A_{K,f}$ of finite adeles of $K$ identifies with
$K \otimes_\Q \A_{\Q,f}$ with $\A_{\Q,f} = \hat \Z \otimes \Q$. We
will write $\A$, resp. $\A_f$, for $\A_{\Q}$, resp. $\A_{\Q,f}$.

Let $\h$ denote the {\it upper half plane} in $\C$, and let
$\h_\pm = \C - \R$. Put
$$
G : \, = \, R_{K/\Q} {\rm GL}(2)/K,
$$
where $R_{K/\Q}$ denotes the Weil restriction of scalars from $K$
to $\Q$. Then $G$ is a reductive algebraic group over $\Q$ with
$G(\Q) = {\rm GL}(2,K)$, and $G(\R) = {\rm GL}(2, K_\infty)$ has a
natural action, by fractional linear transformations, on
$$
{\mathcal D} : = \, K\otimes \C - K\otimes \R \, \simeq \,
\h_\pm^d,
$$
Define
$$
h: \C^\ast \, \rightarrow \, G(\R), \, a+ib \to
\delta\left(\begin{pmatrix}
a &b\\
-b &a
\end{pmatrix}\right),
$$
where $\delta$ denotes the diagonal embeddding of GL$(2, \R)$ in
$G(\R)$. Let $K_\infty$ denote the centralizer of $h(\C^\ast)$ in
$G(\R)$, so that we have the identification
$$
{\mathcal D} \, = \, G(\R)/K_\infty .
$$

Let $C$ be an open compact subgroup of $G_f = G(\A_f)$. Denote by
$S_C = S_C(G,h)$ the associated $d$-dimensional {\it Shimura
variety} over $\Q$, which is quasi-projective and attached to the
datum $(G,h; C)$ by the theory of canonical models ([De1]). One
has
$$
S_C(\C) \, = \, G(\Q)\backslash {\mathcal D} \times G_f/C.
$$
We will take $C$ to be small enough so that $G(\Q) \cap C$ has no
elliptic elements, making $S_C$ non-singular.

The standard approximation theorem for $\mathbb G_m$ says that we
can find a finite set of elements $b_1, b_2, \dots, b_{h(C)}$ in
$\A_{K,f}^\ast$ such that
$$
\A_K^\ast \, = \, \cup_{j=1}^{h(C)} K^\ast b_jK_\infty^+ {\rm
det}(C),
$$
where $K_\infty^+$ denotes the totally positive elements in $K
\otimes \R \simeq \R^d$. Combining this with the strong
approximation theorem for SL$(2)$/K and the finiteness of the
class number, and denoting by $G(\R)^+$ the subgroup of $G(\R)$
consisting of totally positive elements, one obtains the following
useful decomposition:
$$
G(\A) \, = \cup_{j=1}^{h(C)} G(\Q)x_jG(\R)^+ C,
$$
with
$$
x_j \, = \, \begin{pmatrix} b_j & 0\\ 0 & 1
\end{pmatrix}
$$
Consequently, we get an identification of quasi-projective,
complex varieties:
$$
S_C(\C) \, = \, \cup_{j=1}^{h(C)} S_{\Gamma_j}(\C),
$$
where for each $j$,
$$
S_{\Gamma_j}(\C) \, = \, \Gamma_j\backslash \h_\pm^d,
$$
with
$$
\Gamma_j \, = \, G(\Q) \cap x_jg(\R)^+Cx_j^{-1}.
$$
Each $\Gamma_j$ is a discrete subgroup of $G(\R)$ and
$\Gamma_j\backslash \h_\pm^d$ is a Hilbert modular $d$-fold in the
classical sense, having a canonical model $S_{\Gamma_j}$ over a
finite abelian extension $k(\Gamma_j)$ of $\Q$. Gal$(\Q^{\rm
ab}/\Q)$ acts continuously, but not transitively, on the group
$\pi_0(S_C)$ of connected components of $S_C(\C)$ ([De1]).

Put
$$
S : = \, \lim\limits_C \, S_C,
$$
which is a pro-variety over $\Q$ admitting a right $G_f$-action.

Let $S_C^\ast$ denote the projective (singular) {\it Baily-Borel
compactification} of $S_C$ over $\Q$. One has
$$
S_C^\ast \, = \, S_C \cup S_C^\infty,
$$
where $S_C^\infty$ is a finite set of {\it cusps}.

Choose and fix a smooth {\it toroidal compactification} $\tilde
S_C$ over $\Q$ ([AMRT]), defined by a rational cone decomposition.
Let $\tilde S_C^\infty$ stand for the inverse image of
$S_C^\infty$ in $\tilde S_C$, which we can (and we will) arrange
to be a divisor with normal crossings. The irreducible components
are smooth rational varieties, and this will be important to us.
One can construct smooth toroidal compactifications $\tilde
S_{\Gamma_j}(\C)$ of the components $S_{\Gamma_j}(\C)$ such that
$$
\tilde S_C(\C) \, = \, \cup_{j=1}^{h(C)} \tilde S_{\Gamma_j}(\C).
$$

\vskip 0.2in

\section{Contribution from infinity}

\bigskip

Let $K$ be a quartic, totally real extension of $\Q$, and $\ell$ a
prime. Fix a compact open subgroup $C$ as in the section above,
together with a smooth toroidal compactification $\tilde S_C$ over
$\Q$ of the associated Shimura variety $S_C$. Since we are only
interested in the cycles of codimension $2$ on these fourfolds, we
will write $r_{\rm alg}$ for $r_{\rm alg}^2$ and $r_{\rm an}$ for
$r_{\rm an}^2$.

By the {\it decomposition theorem} of Beilinson, Bernstein and
Deligne ([BBD], section 5.4), we have a short exact sequence of
$\mathcal G_\Q$-modules
$$
0 \, \to \, IH^4_{et}(\tilde S_C \times_\Q \overline \Q, \Q_\ell)
\, \to \, H^4_{et}(S_C^\ast \times_\Q \overline \Q, \Q_\ell) \,
\stackrel{b}{\longrightarrow} \, H^4_{\tilde S_C^\infty,
et}(\tilde S_C \times_\Q \overline \Q, \Q_\ell) \, \to \, 0,
\leqno(3.1)
$$
where $IH^\ast$ denotes the middle intersection cohomology of
Goresky, MacPherson and Deligne of the Baily-Borel
compactification $S_C^\ast$, which is pure by a theorem of Gabber;
in other words, the eigenvalues of Frobenius elements $Fr_P$
acting on $IH^4_{et}(S_C^\ast \times_\Q \overline \Q, \Q_\ell)$
have, at good primes $P$, absolute value $N(P)^2$. And the group
on the right of (3.1) signifies the cohomology of $\tilde S_C
\times_\Q \overline \Q$ with supports in the divisor $\tilde
S_C^\infty$.

Though we do not need it here, one knows enough to show explicitly
that the sequence (3.1) splits as ${\mathcal G}_\Q$-modules. We
have recently learnt that a very general result of this sort for
all Shimura varieties has been established by A.~Nair ([N]).

\medskip

Define, for $\nu = $alg, an or $\ell$, and for any number field
$k$, the corresponding ranks $r_{\nu, k}(S_C^\ast)$ and $r_{\nu,
k}(\tilde S_C^\infty)$ by using the respective ${\mathcal
G}_\Q$-modules $IH^4_{et}(\tilde S_C \times_\Q \overline \Q,
\Q_\ell)$ and $H^4_{\tilde S_C^\infty, et}(\tilde S_C \times_\Q
\overline \Q, \Q_\ell)$.

\medskip

\noindent{\bf Proposition 3.2} \, \it We have, for any number
field $k$,
$$
r_{{\rm alg}, k}(\tilde S_C^\infty) \, = \, r_{\ell, k}(\tilde
S_C^\infty)\leqno(a)
$$
and
$$ r_{{\rm an}, k}(\tilde S_C^\infty) \, = \, r_{\ell, k}(\tilde
S_C^\infty)\leqno(b)
$$
\rm

\medskip

{\it Proof of Proposition 3.2}. \,

(a) It suffices to prove this over a sufficiently large extension
of $k$, which contains in particular the (abelian) fields of
definition of the cusps. Since the divisor $\tilde S_C^\infty$ is
the inverse image of the set $S_C^\infty$, and since the cusps are
isolated, we have in \'etale as  well as singular cohomology,
$$
H^\ast_{\tilde S_C^\infty}(\tilde S_C) \, \simeq \, \oplus_{\sigma
\in S_C^\infty} H^\ast_{D_\sigma}(\tilde S_C), \leqno(3.3)
$$
where each $D_\sigma$, the fiber over $\sigma$, is a {\it divisor
with normal crossings}. Let $\{D^i_\sigma \, \vert \, 1 \leq i
\leq r(\sigma)\}$ denote the set of irreducible components of
$D_\sigma$. One knows (cf. [AMRT]) that each $D^i_\sigma$ is a
smooth toric threefold, and so by purity the cohomology of $\tilde
Y$ with supports in $D^i_\sigma$ can be expressed in terms of the
cohomology of $D^i_\sigma$.

Clearly, it suffices to prove the assertion for each cusp
separately. From the geometry of $\tilde D_\sigma$ we obtain the
following exact sequence
$$
\to \, \oplus_{i=1}^{r(\sigma)} H^2(D^i_\sigma)(1) \, \to \,
H^4_{D_\sigma}(\tilde S_C)(2) \, \to \, \oplus_{i \ne j}
H^0(D^{i,j}_\sigma) \, \to \, , \leqno(3.4)
$$
where $D^{i,j}_\sigma$ denotes, for all unequal $i,j$ with $1 \leq
i,j \leq r(\sigma)$, the intersection $D^i_\sigma \cap
D^j_\sigma$. It is clear then that the Tate classes of codimension
$2$ on $\tilde Y$ with supports in $D_\sigma$ are generated by the
following:

\noindent{(3.5)}
\begin{enumerate}
\item[(i)] \, Tate classes of codimension $1$ on the components
$D^i_\sigma$, \, {\it and}
\item[(ii)] \, Classes of the intersections $D^{i,j}_\sigma$.
\end{enumerate}
The classes of type (ii) are obviously algebraic. And so are those
of type (i) because the $D^i_\sigma$ are rational varieties; in
fact, for every $(\sigma, i)$, the entire $H^2(D^i_\sigma)(1)$ is
generated by divisor classes. Hence we get (a).

\medskip

(b) \, It is evident from (3.3), (3.4) and (3.5) that the Galois
representation on $H^4_{\tilde S_C^\infty, et}(\tilde S_C
\times_\Q \overline \Q, \Q_\ell)$ is {\it potentially abelian},
and its $L$-function over any number field $k$ is associated to
that of a $\C$-representation of the absolute Weil group $W_k$. In
this case the equality of (b) is well known (cf. [De2]). Done.

\qed

\medskip

There is a natural analog of the decomposition (3.1) for the Betti
cohomology with $\Q$-coefficients, with $IH^4_B(\tilde S_C(\C),
\Q)$ being a pure $\Q$-Hodge structure of weight $4$. So we may
define the {\it Hodge cycle ranks} $r_{Hg}(S_C^\ast)$ and
$r_{Hg}(\tilde S_C^\infty)$ by using $IH^4_B(\tilde S_C(\C), \Q)$
and $H^4_{\tilde S_C^\infty, B}(\tilde S_C(\C), \Q)$ respectively.

Arguing as in the proof of Proposition 3.2, we also get the
following

\medskip

\noindent{\bf Proposition 3.6} \, \it We have
$$
r_{{\rm alg}}(\tilde S_C^\infty) \, = \, r_{Hg}(\tilde
S_C^\infty).
$$ \rm

\medskip

In view of Propositions 3.2 and 3.6, we see that any Tate or Hodge
class of codimension $2$ on $\tilde S_C$ with supports in $\tilde
S_C^\infty$ is algebraic with the predicted order of pole (in the
Tate case). It then follows from the decomposition (3.1) that to
establish the {\it Tate conjectures} over any $k$, resp. the Hodge
conjecture, for $\tilde S_C$, it suffices to prove the following
identities:
$$
r_{{\rm alg}, k}(S_C^\ast) \, = \, r_{{\rm an}, k}(S_C^\ast),
\leqno(3.7-T1)
$$
and
$$
r_{{\rm alg}, k}(S_C^\ast) \, = \, r_{\ell, k}(S_C^\ast),
\leqno(3.7-T2)
$$
resp.
$$
r_{{\rm alg}}(S_C^\ast) \, = \, r_{Hg}(S_C^\ast). \leqno(3.7-{\rm
Hg})
$$

\vskip 0.2in

\section{Decomposition according to Hecke}

\bigskip

Put
$$
V_B : = \, IH^{4}_B(S_C^\ast(\C), \Q)
\leqno(4.1)
$$
and
$$
V_\ell : = \, IH^4_{et}(S_C^\ast\times_\Q \overline \Q, \Q_\ell).
$$

\medskip

Let $\mathcal H_C$ denote the $\Q$-linear Hecke algebra of level
$C$, which acts semisimply on $V_B$ and $V_\ell$. It is generated
by the characteristic functions of the double cosets $CgC$, $g \in
G_f$. For every field extension $E$ of $\Q$, we will set
$$
\mathcal H_{C,E} \, = \, \mathcal H_C \otimes_\Q E.
$$

The elements of $\mathcal H_C$ act as algebraic correspondences of
finite degree on $S_C$. To be precise, any $g$ in $G_f$ acts on
$S$ on the right, but it does not preserve $S_C$. But if we put
$C_g: = C \cap gCg^{-1}$, it is again an open compact subgroup of
$G_f$, and there are two homomorphisms $C_g \to C$ given by the
identity and the conjugation by $g^{-1}$, resulting in a
corresponding pair of maps, denoted $R(1)$ and $R(g)$, from $S_C$
into $S_{C_g}$. This leads to a self-correspondence, called the
{\bf Hecke correspondence}, of $S_C$ given by
$$
T_g(x) \, = \, R(g)(R(1)^{-1}(x)).
$$
It is not hard to see that $T_g$ depends only on the double coset
$CgC$, and this way one gets an isomorphism of the $\Q$-algebra
generated by the $T_g$ with ${\mathcal H}_C$. The Hecke
correspondences also extend, as ramified correspondences, to the
Baily-Borel compactification $S_C^\ast$. (They do not extend to
the toroidal compactification for a fixed rational cone
decomposition, which is the {\it raison d'etre} for the discussion
of the previous section.) Hence they act on the intersection
cohomology of $S_C^\ast$, and it is natural to decompose according
to isotypic subspaces according to the algebra.

It is a basic fact that irreducibles $\eta$ of $\mathcal H_{C,
\overline \Q}$ are in bijection with $\overline \Q$-irreducible
admissible representations $\pi_f$ of $G_f$ which admit a non-zero
$C$-fixed vector, the correspondence being given by $\eta =
\pi_f^C$. Let $E = E(\pi_f)$ denote the field of rationality of
$\pi_f$; it is known that $\pi_f$ and $\eta = \pi_f^C$ can be
realized over $E$. In sum, irreducibles of $\mathcal H_C$ are
parametrized by the packets
$$
\{\pi_f^\sigma \, \vert \, \pi_f^C \ne 0, \, \sigma \in {\rm Hom}(E, \overline
\Q)\}
$$
One gets the decomposition for $\ast = B, et$:
$$
V_\ast \, \simeq \, \oplus_{\pi_f \in \Sigma_C } \,
V_\ast(\pi_f)^{m(\pi_f,C)},
\leqno(4.2)
$$
where $\Sigma_C$ denotes the set, modulo Galois conjugation, of
irreducible admissible $\overline \Q$-representations $\pi_f$ of
$G_f$ admitting a non-zero vector fixed by $C$, and $m(\pi_f, C)$
denotes the dimension of $C$-invariants. There is an isomorphism
of $\Q_\ell$-vector spaces:
$$
V_B(\pi_f) \otimes_\Q \Q_\ell \, \simeq \, V_\ell(\pi_f).
\leqno(4.3--a)
$$
Moreover, each $V_B(\pi_f)$ admits an $E$-action. It is of
importance to know the dimension of $V_B(\pi_f)$ over $E$ and we
will come to this question below. The $\ell$-adic Galois
representation is an $E \otimes \Q_\ell$-module and it decomposes
as
$$
V_\ell(\pi_f) \, \simeq \, \oplus_{\lambda \vert \ell} \,
V_\lambda(\pi_f). \leqno(4.3--b)
$$

One can define in the obvious way, for any number field $k$, the
$\pi_f$-components of the algebraic and analytic ranks $r_{{\rm
alg}, k}(S(\pi_f))$ and $r_{{\rm an}, k}(S(\pi_f))$ respectively.
(We want to think of $S(\pi_f)$ as a Grothendieck motive with
coefficients in $E$.)

\medskip

We will say that $\pi_f$ is $CM$ (or of {\it CM type}) iff $\pi_f
\simeq \pi_f \otimes \nu$, for a non-trivial quadratic character
$\nu$. Then $\pi_f$ is necessarily infinite-dimensional and the
quadratic extension $M = M(\nu)$ of $F$ cut out by $\nu$ is a CM
field. Put, for $\ast = B, et$,
$$
V_\ast^{CM} \, = \, \oplus_{\pi_f \in \Sigma_C, \, \pi_f \, CM} \,
V_\ast(\pi_f)^{m(\pi_f,C)}, \leqno(4.4--a)
$$
and
$$
V_\ast' \, = \, V_\ast/V_\ast^{CM}, \leqno(4.4--b)
$$

\medskip

In view of Proposition 3.2 and the decomposition (4.2) above,
Theorem A will be a consequence of the following:

\medskip

\noindent{\bf Theorem A$^\prime$} \, \it Suppose $K$ is normal
over $\Q$, and $k$ an abelian number field. Let $\pi_f$ be an
irreducible, admissible, $\overline \Q$-rational, non-CM
representation of $G(\A_f)$, equipped with a non-zero vector in
(the space of) $\pi_f$ fixed by the compact open subgroup $C$,
such that $V_B(\pi_f) \ne 0$. Then we have
$$
r_{{\rm alg}, k}(S(\pi_f)) \, = \, r_{{\rm an}, k}(S(\pi_f)).
$$
\rm

\medskip

Put
$$
V_{(2)} \, = \, H^4_{(2)}(S_C(\C), \C),
$$
where the group on the right is the {\it $L^2$-cohomology} of the
open manifold $S_C(\C)$. One knows by the proof of Zucker's
conjecture by Saper-Stern ([SaSt] and Looijenga ([Lo]), that there
is an isomorphism
$$
V_B \otimes_\Q \C \, \simeq V_{(2)}.
\leqno(4.5)
$$

If $\mathfrak G_\C$ denotes the set of matrices in the
complexified Lie algebra of $G$ with purely imaginary trace, there
is a standard isomorphism ([BoW])
$$
H^\ast_{(2)}(S_C(\C), \C) \, \simeq \, H^\ast({\mathfrak G}_\C,
K_\infty; L^2_{\rm disc}(G(\Q)Z(\R)\backslash G(\A)/C)^\infty),
\leqno(4.6)
$$
where $L^2_{\rm disc}(G(\Q)Z(\R)\backslash G(\A)/C)$ denotes the
$C$-invariants in the discrete spectrum of (the right regular
representation) $L^2(G(\Q)Z(\R)\backslash G(\A))$, and $Z$ the
center of $G$. The superscript $\infty$ signifies taking the
subspace of smooth vectors at infinity. One has a decomposition as
unitary $G(\A)$-modules
$$
L^2_{\rm disc}(G(\Q)Z(\R)\backslash G(\A)) \, \simeq \, L^2_{\rm
res}(G(\Q)Z(\R)\backslash G(\A)) \hat\oplus L^2_{\rm
cusp}(G(\Q)Z(\R)\backslash G(\A)), \leqno(4.7)
$$
where the second space on the right is the {\it space of cusp
forms} ([BoJ]), while the  first space on the right is spanned by
the {residual representations}, which in this case are precisely
the one-dimensional unitary representations $\pi$ occurring in the
discrete spectrum. And $\hat{\oplus}$ signifies taking the Hilbert
direct sum.

Using (4.5), (4.6, the complete reducibility of $L^2_{\rm
disc}(G(\Q)Z(\R)\backslash G(\A))$, and the fact that this
representation is {\it multiplicity free}, we get
$$
V_B \otimes_\Q \C \, \simeq \, \hat\oplus_{\pi} \, (H^4({\mathfrak
G}_\C, K_\infty; {\mathcal H}(\pi_\infty)) \otimes \pi_f^C),
\leqno(4.8)
$$
where $\pi$ runs over the irreducible unitary representations of
$G(\A)$ (up to equivalence) (admitting a $C$-fixed vector), and
${\mathcal H}_\pi$ denotes the space of smooth vectors of
$\pi_\infty$.

Let Coh$_{G(\R)}$ denote the set of all irreducible unitary
representations of $G(\R)$ (up to equivalence) with trivial
central character and with non-zero $({\mathfrak G}_\C,
K)$-cohomology in degree $4$, and let Coh$_{G}(C)$ denote the set
consisting of the $\pi$ occurring in $L^2_{\rm
disc}(G(\Q)Z(\R)\backslash G(\A)$ such that $\pi_\infty \in
$Coh$_{G(\R)}$ and $\pi_f^C \ne 0$.

By the strong multiplicity one theorem, any cuspidal automorphic
representation $\pi$ is determined by the knowledge of almost all
of its local components, in particular by its finite part $\pi_f$.
The analogous statement about the one-dimensional $\pi$ is
obvious. In view of this one gets, by comparing (4.2) and (4.8)
for each $\pi \in $Coh$_G(C)$,
$$
(V_B(\pi_f) \otimes_\Q \C)^{{\rm dim} \pi_f^C} \, \simeq \,
H^4({\mathfrak G}_\C, K_\infty; {\mathcal H}(\pi_\infty) \otimes
\pi_f^C. \leqno(4.9)
$$

The Galois conjugate of any cuspidal , resp. one dimensional,
$\pi$ in Coh$_G(C)$ is again cuspidal, resp. one dimensional. For
$\ast = B$ or $et$, let us set
$$
V_\ast^{\rm res} \, \simeq \, \oplus_{\pi \in {\rm Coh}_G(C), {\rm
dim}(\pi) = 1} \, V_\ast(\pi_f)^{m(\pi_f,C)},
\leqno(4.10)
$$
and
$$
V_\ast^{\rm cusp} \, \simeq \, \oplus_{\pi\in {\rm Coh}_G(C), \pi
\, {\rm cuspidal}} \, V_\ast(\pi_f)^{m(\pi_f,C)}.
$$

\medskip

\noindent{\bf Proposition 4.11} \, \it Let $\pi \in $Coh$_G(C)$
with dim$(\pi) = 1$. Then the Hodge and Tate classes in $V^{\rm
res}(2)$ are algebraic, in fact represented by Chern classes of
vector bundles. And we have, for any number field $k$,
$$
r_{\rm alg}(S(\pi_f)) \, = \, r_{\rm an}(S(\pi_f)).
$$
\rm

\medskip

{\it Proof}. \, Recall from section 2 that $S_C$ is a finite union
of connected Hilbert modular varieties $S_{\Gamma}$, whose complex
points are quotients of $\mathcal H_\pm^4$ by a congruence
subgroup $\Gamma$ of $G(\R)$. Restriction of any cohomology class
in $V_B^{\rm res}$ to $S_\Gamma$ comes from the continuous
cohomology of $G(\R)$, i.e., represented by a $G(\R)$-invariant
differential form on $\mathcal H_\pm^4$. If $z =
(z_1,z_2,z_3,z_4)$, $z_j = x_j +iy_j$, represents a point on
$\mathcal H_\pm^4$, then any such $G(\R)$-invariant differential
$4$-form is spanned by the following forms (with $1 \leq j,r \leq
4$, $j \ne r$):
$$
\omega_{j,r} \, = \, \eta_j \wedge \eta_r
$$
where
$$
\eta_j \, = \, \frac{dz_j{\overline {dz}}_j}{y_j^2}.
$$
Each $\eta_j$ is just the volume form on the $j$-th factor. It is
well known that $\eta_j$ is the Chern class of a line bundle,
namely the one given by the tangent bundle on $\mathcal H_\pm$. It
represents a divisor on $\mathcal H_\pm^4$, which descends to one
on $S_\Gamma$; call it $D_j$. The intersection of $D_j$ and $D_r$
gives, for $j \ne r$, a codimension $2$ cycle $Z_{j,r}$
represented by a class in $V_B^{\rm res}$. (To be precise, the
Chern class lies in $V_B^{\rm res}(2)$.) The associated Galois
representation on $V_\ell^{\rm res}$ is potentially abelian, and
it is easy to match the poles of this $L$-function (for any $k$)
with the Tate classes, i.,e., the Galois invariants in
$V_{\ell}^{\rm res}(2)$. Moreover, these Tate cycles, just like
the Hodge cycles, are represented by the Chern classes. The
Proposition follows.

\qed

\medskip

So, in order to prove Theorem $A'$ (and hence Theorem A), we may,
and we will, concentrate on the {\it cuspidal case} from here on.

\medskip

For every ideal $\mathfrak N$ of $K$ with prime factorization
$\prod_v {\mathfrak P}_v^{f_v}$, denote by $C_1({\mathfrak N})$
the compact open subgroup of $G_f$ given by
$$
C_1({\mathfrak N}) \, = \, \prod_v \, C_1({\mathfrak P}_v^{f_v})
\leqno(4.12)
$$
with
$$
C_1({\mathfrak P}_v^{f_v}) \, = \, \{\begin{pmatrix} a_v & b_v\\
c_v & d_v\end{pmatrix} \, \vert \, c_v, d_v-1 \in {\mathfrak
P}_v^{f_v}\}.
$$
One knows that given any cuspidal automorphic representation
$\pi$, it has a {\it conductor} ${\mathfrak N}(\pi)$, which is the
largest ideal ${\mathfrak N}$ such that the space of $\pi_f$ has a
non-zero vector fixed by $C_1({\mathfrak N})$. One also knows that
the space of ${\mathfrak N}(\pi)$-invariants in $\pi_f$ is exactly
one dimensional. So we get, for any cuspidal $\pi$ in Coh$_G(C)$:
$$
V_B(\pi_f) \otimes_\Q \C \, \simeq \, H^4({\mathfrak G}_\C,
K_\infty; {\mathcal H}(\pi_\infty)) \otimes \pi_f^{C(\pi)}.
\leqno(4.13)
$$
where $C(\pi)$ is shorthand for $C_1({\mathfrak N}(\pi))$.

It is known that given any cuspidal $\pi$ contributing to the
cohomology in degree $4$, we must have
$$
\pi_\infty \, \simeq \, {\mathcal D}_2^{\otimes 4},
\leqno(4.14)
$$
where ${\mathcal D}_2$ is the discrete series representation of
GL$(2, \R)$ of lowest weight with trivial central character, which
contributes to cohomology in degree $1$ only. Consequently,
$$
H^4({\mathfrak G}_\C, K_\infty; {\mathcal H}(\pi_\infty)) \,
\simeq \, H^1({\mathfrak Gl}(2, \C), L_\infty; {\mathcal
D}_2^\infty)^{\otimes 4},
\leqno(4.15)
$$
where $L_\infty \simeq $ SO$(2)\R_+^\ast$. It is well known that
$H^1({\mathfrak Gl}(2, \C), L_\infty; {\mathcal D}_2^\infty)$ is
$2$-dimensional. Hence it follows that
$$
{\rm dim}_E V_B(\pi_f) \, = \, 16.
\leqno(4.16)
$$
If $V_B(\pi_f)_\C^{p,q}$ denotes the Hodge $(p,q)$-piece of
$V_B(\pi_f) \otimes \C$, then we get
$$
{\rm rk}_{E \otimes \C} V_B(\pi_f)_\C^{2,2} \, = \, 6.
\leqno(4.17)
$$
Applying Propositions 1.19 and 1.20, we then get (for any number
field $k$)
$$
r_{{\rm alg}, k}(\pi_f) \, \leq \, r_{\ell, k}(\pi_f) \, \leq \,
6. \leqno(4.18)
$$
We will see later that for $k$ abelian, and $\pi$ cuspidal and
non-CM), $r_{\ell, k}(\pi_f)$ is at most $2$. For $\pi$ of CM
type, however, it could be more than $2$ for $k$ abelian, and
could be $6$ for suitable (non-abelian) $k$ when $K$ is
biquadratic.

\vskip 0.2in

\section{Twisted analogues of Hirzebruch-Zagier cycles}

\bigskip

In this section $K$ will denote a totally real number field of
degree $m$ and $F$ a subfield. Put $r = [K:F]$. Put
$$
G \, = \, R_{K/\Q} {\rm GL}(2)/K \quad {\rm and} \quad H \, = \,
R_{F/\Q} {\rm GL}(2)/F. \leqno(5.1)
$$
The map $h$ defined in section 2 factors through a map $h_0$ of
$R_{\C/\R}(\C^\ast)$ into $H_\R$. Let $S$, resp. ${}^HS$, denote
the Shimura variety over $\Q$ associated to $(G, h)$, resp. $(H,
h_0)$, with Baily-Borel compactifications $S^\ast, {}^HS^{\ast}$
respectively. The natural embedding of $H$ into $G$ leads to an
embedding (over $\Q$)
$$
\delta: {}^HS \, \hookrightarrow \, S,
$$
which extends to a map from ${}^HS^\ast$ into $S^\ast$.

Recall that $S$ comes equipped with a $G_f$-action by right
translation $R$. For any open compact subgroup $C$ of $G_f$, let
$$
p_C: \, S \, \rightarrow \, S_C
$$
denote the natural projection, which extends to $S^\ast
\rightarrow S_C^\ast$. For any $g \in G_f$, define the
corresponding {\bf Hirzebruch-Zagier cycle}, or {\bf HZ-cycle} for
short, (relative to $H$) to be the {\it algebraic cycle of
codimension $m-r$} given by
$$
{}^FZ_{g,C} \, = \, p_C(R(g)(\delta({}^HS))) \leqno(5.2)
$$
with compactification
$$
{}^FZ_{g,C}^\ast \, = \, p_C(R(g)(\delta({}^HS^\ast))).
$$

Note that if $x, y \in {}^HS^\ast$ are in the same orbit under
$H_f \cap gCg^{-1}$, which is a compact open subgroup of $H_f$,
then they have the same image in ${}^FZ_{g,C}^\ast$. Thus one
obtains a non-trivial morphism over $\Q$:
$$
{}^HS^\ast_{H_f \cap gCg^{-1}} \, \rightarrow \, {}^FZ_{g,
C}^\ast. \leqno(5.3)
$$

The right translation action of $g$ on $S$ does not descend to a
morphism $S_C \to S_C$, but it does define the Hecke
correspondence $T_g$ encountered earlier. We can view ${}^FZ_{g,
C}$ as the image of the Hilbert modular subvariety ${}^HS_{C \cap
H_f}$ under $T_g$.

These cycles, or rather their classical versions of them, were
introduced by Hirzebruch and Zagier ([HZ]) in the case $F=\Q$ and
$K$ real quadratic, which were used in [HLR] to prove the Tate
conjecture for Hilbert modular surfaces over abelian fields.

\medskip

Now on to the {\it twisted versions}. Let $\mu$ be any finite
order character of (the idele class group of) $K$ of conductor
$\mathfrak c$. Let $C = C_1(\mathfrak N)$. Put
$$
C[\mu] \, = \, C_1({\rm lcm}(\mathfrak N, \mathfrak c^2))
\leqno(5.4)
$$
(see (4.12)). Write $\mathfrak O_{\mathfrak c}$ for the ring of
integers of $F_{\mathfrak c}: = \prod_{v \mid \mathfrak c} F_v$,
and define a subset of $F_{\mathfrak c}$ by
$$
X: = \, \{x = (x_v) \in F_{\mathfrak c} \, \vert \, v(x_v) \geq
-v(\mathfrak c), \forall \, v\}. \leqno(5.5)
$$
Let $\tilde X$ be a set of representatives in $X$ for $X$ mod
$\mathfrak O_{\mathfrak c}$, which is a group isomorphic to
$\mathfrak O_{\mathfrak c}/\mathfrak c\mathfrak O_{\mathfrak c}$.
To each $t$ in $\tilde X$, associate the unipotent matrix $u(t) =
\begin{pmatrix} 1 & t\\ 0 & 1 \end{pmatrix}$.
Now recall from section 2 that
$$
S_C(\C) \, = \, \cup_{j=1}^{h(C)} \Gamma_j\backslash \h_\pm^d,
$$
with
$$
\Gamma_j \, = \, G(\Q) \cap x_jg(\R)^+Cx_j^{-1}.
$$

Recall that every $x$ in $G_f$ defines a {\it Hecke
correspondence} $T(x)$ of $S_C$, which does not in general
preserve the connected components $S_C^{(j)}$ of $S_C$. Let
$T_j(x)$ denote the restriction of $T(x)$ to $S_C^{(j)}$ when
det$(x) = 1$.

The {\it twisting correspondence} $R(\mu) \subset S_C \times
S_{C[\mu]}$ can now be defined (cf. [MuRa], section 2, for
example) as
$$
R(\mu) \, = \, \sum_{j=1}^{h(C)} \mu_f({\rm det}(x_j)) R_j(\mu),
\leqno(5.7)
$$
with
$$
R_j(\mu) \, = \, \sum_{t \in \tilde X} \, T_j(u_t).
$$
It is easy to verify that for all $x \in G_f$,
$$
T(x) \circ R(\mu) \, = \, \mu_f({\rm det}(x))R(\mu)\circ T(x).
$$

The twisting correspondence, being algebraic, acts on any
cohomology group, Betti or \'etale, of the fourfold $S = {\rm lim}
S_C$. The induced operator sends the $\pi_f$-component to the
$\pi_f \otimes \mu$-component. It may be useful to note that the
twisting correspondence $R(\mu)$ is rational over $\Q(\mu_1)$,
where $\mu_1$ denotes the restriction of $\mu$ to (the ideles of)
$\Q$.

\medskip

Given a H-Z cycle $Z$ on $S$ and a character $\mu$, define the
associated {\it $\mu$-twisted H-Z cycle} $Z(\mu)$ to be the
push-forward of $Z$ under  $R(\mu)$. It is again algebraic and
rational over $\Q(\mu_1)$.

\vskip 0.2in

\section{Asai ${\bf L}$-functions of type $(n,d)$}

\bigskip

Fix $n \geq 1$ and let $K/F$ be an extension of number fields of
degree $d$ with Galois closure $\tilde K$ (over $F$). Then
Gal$(\tilde K/F)$ acts on ${\rm Hom}(K,\C)$ and hence on the
$\C$-vector space $\C^{{\rm Hom}(K,\C)} \simeq \C^d$. This in turn
induces an action on the group ${\rm GL}(n,\C)^d$.

Now consider the algebraic $F$-group
$$
G \, = \, R_{K/F} {\rm GL}(n)/K, \leqno(6.1)
$$
where $R_{K/F}$ is the Weil restriction of scalars, so that $G(F)
= {\rm GL}(n,K)$. Its {\it dual group} can be taken to be the
semidirect product
$$
{}^LG \, = \, {\rm GL}(n,\C)^d \bowtie W_F, \leqno(6.2)
$$
where the absolute Weil group $W_F$ acts via its quotient
$W_F/W_{\tilde K} \simeq {\rm Gal}(\tilde K/F)$, the action being
the one described above. Then we have the natural representation
$$
r_{K/F}: \, {}^LG \, \rightarrow \, {\rm GL}((\C^n)^{\otimes d})
\, \simeq \, {\rm GL}(n^d, \C), \leqno(6.3)
$$
given by
$$
r_{K/F}((g_\sigma); 1)(\otimes_{\sigma} v_\sigma) \, = \,
\otimes_\sigma  g_\sigma v_\sigma \leqno(6.4)
$$
and
$$
r_{K/F}((e_\sigma) ; \tau)(\otimes_\sigma v_\sigma) \, = \,
\otimes_\sigma  v_{\tau\sigma},
$$
for all $(g_\sigma)$ in GL$(n,\C)^{{\rm Hom}(K,\C)}$ and $\tau \in
W_F$, where $\sigma$ runs over Hom$(K, \C)$ and $e_\sigma$ denotes
the identity $n \times n$-matrix matrix in the $\sigma$-th place.

For any automorphic representation $\pi = \otimes'_v \pi_v$ of
GL$(n, \A_K)$, and for any idele character $\chi$ of $F$, which we
may view by class field theory as a character, again denoted
$\chi$, of $W_F$ (and hence of ${}^LG$), let $L(s, \pi; r_{K/F}
\otimes \chi)$ denote the associated Langlands $L$-function, which
we will call the {\it $\chi$-twisted Asai $L$-function of $\pi$
relative to $K/F$}, which is of degree $n^d$ over $F$. It has an
Euler product expansion over the places $v$ of $F$, convergent in
a right half plane, with local factors $L_v(s, \pi; r_{K/F}\otimes
\chi )$ and $\varepsilon_v(s, \pi; r_{K/F}\otimes \chi)$. If $v$
is a finite place of $F$ such that $\pi_w$ is unramified at any
place $w$ of $K$ above $v$, then there is a semisimple conjugacy
class $A(\pi_v)$ in ${}^LG$ such that
$$
L(s, \pi_v; r_{K/F}\otimes \chi_v) \, = \, {\rm
det}(I-\chi_v(\varpi_v)r_{K/F}(A(\pi_v))Nv^{-s})^{-1},
\leqno(6.5)
$$
where $\chi_v$ denotes the $v$-component of $\chi$, which is a
character of $W_{F_v} \simeq F_v^\ast$, and $\varpi_v$ denotes the
uniformizer at $v$. Clearly, this $L$-factor is $1$ unless
$\chi_v$ is also unramified.

In order to describe the local factors at all the places, which
was originally done by Langlands at the places which are
unramified for the datum, we need some preliminaries. At any place
$w$ of $K$, let $W_{K_w}'$ denote the Weil group $W_{K_w}$ if $w$
is archimedean and $W'_{K_w} \times {\rm SL}(2, \C)$ if $w$ is
non-archimedean. One knows by the {\it local Langlands
conjecture}, established long ago over archimedean fields by
Langlands ([La3]), and recently proved for GL$(n)$ over $p$-adic
fields in the independent works of M.~Harris and R.L.~Taylor
([Ha-T]) and Henniart ([He]), that $\pi_w$ is associated to an
$n$-dimensional $\C$-representation $\sigma_w$ of $W'_{K_w}$. This
association $\pi_w \to \sigma_w$ is functorial for taking
contragredients, pairing the central character $\omega_w$ of
$\pi_w$ with the determinant of $\sigma_w$, such that
$$
L(s, \pi_w \times \pi'_w) \, = \, L(s, \sigma_w \otimes \sigma'_w)
\leqno(6.6)
$$
and
$$
\varepsilon(s, \pi_w \times \pi'_w) \, = \, \varepsilon(s,
\sigma_w \otimes \sigma'_w),
$$
for all irreducible admissible representations $\pi'_w$ of GL$(m,
K_w)$, and for all $m$-dimensional representations $\sigma'_w$, $m
\leq n-1$, of $W'_{K_w}$.

Now let $v$ be a place of $F$ and $w$ a place of $K$ above it. Let
$d(w/v)$ denote $[K_w:F_v]$, so that $d = \sum\limits_{w \vert v}
\, d(w/v)$. Let $M_{K_w}^{F_v}(\sigma_w)$ denote the {\bf tensor
induction} (or {\it multiplicative induction}) of $\sigma_w$ from
$W'_{K_w}$ to $W'_{F_v}$ (see [Cu-R], and also [Mu-P]). It is an
$n^{d(w/v)}$-dimensional representation of $W'_{F_v}$, with the
following property:
$$
{\rm Res}_{\tilde K_w}^{F_v}(M_{K_w}^{F_v}(\sigma_w)) \, \simeq \,
\otimes_{\tau \in {\rm Hom}(K_w, \overline K_w)} \sigma_w^\tau,
\leqno(6.7)
$$
where $\tilde K_w$ denotes the Galois closure of $K_w$ over $F_v$.
It is not hard to see that the tensor representation on the right
of (6.7) extends non-uniquely to a representation, and the key
point is that the tensor induction $M_{K_w}^{F_v}(\sigma_w)$ is a
{\it canonical extension}.

Now put
$$
As_{K/F}(\sigma)_v \, = \, \otimes_{w \vert v}
M_{K_w}^{F_v}(\sigma_w), \leqno(6.8)
$$
which is an $n^d$-dimensional representation of $W'_{F_v}$,
$$
L_v(s, \pi; r_{K/F}\otimes \chi ) \, = \,
L(s,As_{K/F}(\sigma)_v)\otimes \chi_v) \leqno(6.9)
$$
and
$$
\varepsilon_v(s, \pi; r_{K/F}\otimes \chi ) \, = \,
\varepsilon(s,As_{K/F}(\sigma)_v \otimes \chi_v)).
$$
Taking the Euler product of (6.9) over all $v$, we get the
definition of the global $\chi$-twisted Asai $L$-function $L(s,
\pi; r_{K/F}\otimes \chi )$.

\medskip

Now suppose $\chi$ is the restriction of an idele class character
$\tilde \chi$ of $K$. (On the Galois side, this corresponds to
{\it transfer}.) Then we have
$$
L(s, \pi \otimes \tilde \chi; r_{K/F}) \, = \, L(s, \pi; r_{K/F}
\otimes \chi). \leqno(6.10)
$$

\medskip

One knows nothing in general about the expected properties of
these $L$-functions (for $d > 1$) except when $d=2$. Here $K/F$ is
a quadratic extension with non-trivial automorphism $\theta$, and
if $\mu$ is a unitary character of $F$, then a cuspidal
automorphic representation $\pi$ of GL$(2, \A_K)$ is said to be
{\it $\mu$-distinguished} ([HLR]) iff the following {\it
$\mu$-period integral} is non-zero for some function $f$ in
${\mathcal V}_\pi$:
$$
{\mathcal P}_\mu(f) :\, = \, \int_{H(F)Z_H(F_\infty)^+\backslash
H(\A_F)} \mu({\text det}(h)f(h) dh,
$$
where $H$ denotes GL$(2)/F$ with center $Z_H$, and $dh$ is the
quotient measure induced by the Haar measure on $H(\A_F)$. When
$\mu = 1$, we will simply say {\it distinguished} when such a
non-vanishing occurs. One has the following well known result:

\medskip

\noindent{\bf Theorem 6.11} \, \it Let $K/F$ be a quadratic
extension of number fields, $n$ a positive integer, and $\pi$ a
cuspidal automorphic representation of GL$(n, \A_F)$ with
contragredient $\pi^\vee$. Then
\begin{enumerate}
\item[(a)] $L(s, \pi; r_{K/F})$ admits a meromorphic continuation
to the whole $s$-plane with a functional equation of the form
$$
L(1-s, \pi^\vee; r_{K/F}) \, = \, \varepsilon(s, \pi; r_{K/F})
L(s, \pi; r_{K/F}),
$$
where $\varepsilon(\pi; r_{K/F})$ is an invertible holomorphic
function.
\item[(b)] \, If $S$ is a finite set of places of $F$ containing
the archimedean places and the primes where $\pi$ is ramified,
$L^S(s, \pi; r_{K/F})$ has a pole at $s=1$ iff $\pi$ is
distinguished.
\end{enumerate}
\rm

\medskip

Part (a) follows from the work of Shahidi ([Sh]) via the
Langlands-Shahidi method, and part (b) was shown by Flicker in
[F$\ell$] by adapting the Rankin-Selberg method of Jacquet,
Piatetski-Shapiro and Shalika ([JPSS]). Note that when $L^S(s,
\pi; r_{K/F})$ has a pole at $s=1$, we have $\pi^\vee \simeq
\pi\circ\theta$.

\vskip 0.2in

\section{Zeroing in on the ${\bf (2,4)}$-case}

\bigskip

The main result of this section is the following:

\medskip

\noindent{\bf Theorem 7.1} \, \it Let $K/F$ be a quartic extension
of number fields such that there is an intermediate field $E$ with
$[K:E] = [E:F] = 2$. Suppose $\pi$ is a cuspidal automorphic
representation of GL$(2, \A_K)$. Then $L(s, \pi; r_{K/F})$ admits
a meromorphic continuation to the whole $s$-plane with a
functional equation of the form
$$
L(1-s, \pi^\vee; r_{K/F}) \, = \, \varepsilon(s, \pi; r_{K/F})
L(s, \pi; r_{K/F}),
$$
where $\varepsilon(\pi; r_{K/F})$ is an invertible holomorphic
function. \rm

\medskip

Thanks to the identity (6.10), we also get the meromorphic
continuation and functional equation of $L(s, \pi; r_{K/F} \otimes
\nu)$ for any idele class character $\nu$ of $F$.

\medskip

When $K/F$ is Galois and $\pi$ of trivial central character and
square-free conductor, it can be shown (see Proposition 8.22 where
$F = \Q$) that $L^S(s, \pi; r_{K/F}\otimes \nu)$ admits, for $S$ a
finite set of places of $F$ containing the archimedean places and
the primes where $\pi$ is ramified, a pole at $s=1$ iff a suitable
twist of $\pi$ is a base change from a quadratic subfield $L$
containing $F$. The order of pole is $2$ iff a twist of $\pi$ is a
base change from $F$ {\it and} $K/F$ is biquadratic.

\medskip

An important ingredient of proof is the following basic result for
quadratic extensions, which we established in [Ra3]:

\medskip

\noindent{\bf Theorem 7.2} \, \it Let $K/E$ be a quadratic
extension of number fields, and let $\pi$ be a cuspidal
automorphic representation of GL$(2, \A_K)$. Then there is an
isobaric automorphic representation $As_{K/E}(\pi)$ of GL$(4,
\A_E)$ such that
$$
L(s, As_{K/E}(\pi)) \, = \, L(s, \pi; r_{K/E}).
$$
\rm

\medskip

One of the steps of the proof of this Theorem in [Ra3] was the
integral representation for this Asai $L$-function when twisted by
a cusp form on GL$(2)/E$. Recently, a different proof has been
given in [Kr] using instead the Langlands-Shahidi theory of this
$L$-function.

\medskip

{\it Proof of Theorem 7.1}. \, The simple reason is that the Asai
$L$-functions can be built in stages. To be precise, we have the
following:

\medskip

\noindent{\bf Proposition 7.3} \, \it Let $K \supset E \supset F$
and $\pi$ be as in Theorem 7.1, and let $As_{K/E}(\pi)$ be as in
Theorem 7.2. Then we have
$$
L(s, \pi; r_{K/F}) \, = \, L(s, As_{K/E}(\pi); r_{E/F}).
$$
\rm

\medskip

{\it Proof of Proposition}. \, It suffices to prove the equality
of local factors. Fix any place $v$ of $F$. Let $u$ be a place of
$E$  above $v$, and $w$ a place of $K$ above $u$. Denote by
$\sigma_w$ the $2$-dimensional representation of $W'_{K_w}$
associated to $\pi_w$ by the local correspondence. Then by the
fact that tensor induction can be achieved in steps (cf, [CR]), we
get
$$
M_{K_w}^{F_v}(\sigma_w) \, \simeq \,
M_{E_u}^{F_v}(M_{K_w}^{E_u}(\sigma_w)). \leqno(7.4)
$$
In view of (6.9), this proves immediately the assertion when $v$
is inert in $K$. Suppose that $v$ is inert or ramified in $E$ with
unique divisor $u$ there, but that $u$ splits into $w, w'$ in $K$.
Then $K_w \simeq K_{w'} \simeq E_u$, and by definition (see (6.8),
(6.9)),
$$
L(s, \pi; r_{K/F}) \, = \, L(s, M_{K_w}^{F_v}(\sigma_w) \otimes
M_{K_{w'}}^{F_v}(\sigma_w)), \leqno(7.5)
$$
while by [Ra3],
$$
As_{K/E}(\pi_w) \, \simeq \, \pi_w \boxtimes \pi_{w'},
$$
where $\boxtimes$ denotes the automorphic tensor product on GL$(2)
\times $GL$(2)$, constructed in [Ra2]. Since $\pi_w \boxtimes
\pi_{w'}$ corresponds to $\sigma_w \otimes \sigma_{w'}$ on the
Weil group side, we obtain
$$
L(s, As_{K/E}(\pi); r_{E/F}) \, = \, L(s, M_{E_u}^{F_v}(\sigma_w
\otimes \sigma_{w'})). \leqno(7.6)
$$
Noting the isomorphism
$$
M_{K_w}^{F_v}(\sigma_w) \otimes M_{K_{w'}}^{F_v}(\sigma_w) \,
\simeq M_{E_u}^{F_v}(\sigma_w \otimes \sigma_{w'}), \leqno(7.7)
$$
we get the assertion of the Proposition when $v$ has a unique
divisor $u$ in $E$, but $u$ splits in $K$. The remaining cases are
similar and are left to the reader.


\vskip 0.2in

\section{Tate classes and the inequality ${\bf r_{\rm alg} \leq r_{\rm an}}$}

\bigskip

Let $K$ be a totally real, quartic extension of $\Q$ containing a
quadratic subfield $F$, and $C$ a compact open subgroup of $G_f$.
Suppose $\pi = \pi_\infty \otimes \pi_f$ is a cuspidal automorphic
representation of $G(\A)$ contributing to Coh$_G(C)$.

\medskip

\noindent{\bf Proposition 8.1} \, \it We have the inequality
$$
r_{{\rm alg}, k}(S(\pi_f)) \, \leq \, r_{{\rm an}, k}(S(\pi_f)),
$$
for any abelian number field $k$. \rm

\medskip

In view of Proposition 1.18, this Proposition will be proved once
we establish the following

\medskip

\noindent{\bf Proposition 8.2} \, \it We have, for any abelian
number field $k$,
$$
r_{\ell, k}(S(\pi_f)) \, = \, r_{{\rm an}, k}(S(\pi_f)).
$$
\rm

\medskip

Before beginning the proof, let us note the following fact which
will be used later:

\medskip

\noindent{\bf Lemma 8.3} \, \it Let $\pi = \pi_\infty \otimes
\pi_f$ be a cusp form on GL$(2)/K$ such that $\pi_f$ admits a
non-zero vector fixed by $C_0({\mathfrak N})$, with ${\mathfrak
N}$ square-free.
\begin{enumerate}
\item[(i)] \, If $P$ is a prime dividing the conductor ${\mathfrak
c}(\pi)$ of $\pi$, then the local component $\pi_P$ must be an
unramified twist of the Steinberg representation of GL$(2, F_P)$.
\item[(ii)] \, If $\pi$ is moreover of CM type, it must be
unramified at every finite place.
\end{enumerate} \rm

\medskip

{\it Proof}. \, By assumption, ${\mathfrak c}(\pi)$ divides
${\mathfrak N}$; then so does the conductor ${\mathfrak
c}(\omega)$ of the central character $\omega$. Suppose $\omega$ is
ramified. Then if $x$ is a new vector in the space of $\pi_f$, the
group $C_0({\mathfrak c}(\omega))$ will act on $x$ by a {\it
non-trivial} character determined by $\omega$, which is trivial on
$C_1({\mathfrak c}(\omega))$. Hence $C_0({\mathfrak N})$ cannot
act trivially on any non-zero vector in $\pi_f$, contradiction! So
$\omega$ must be unramified.

Suppose $P$ is a prime divisor of ${\mathfrak c}(\pi)$. If $\pi_P$
were supercuspidal, then $P^2$ would divide the conductor (see
[Ge], p.73), so the square-freeness assumption forces $\pi_P$ to
be special or a ramified principal series. In the latter case, it
is defined by two local characters $\mu, \nu$ such that $\mu\nu =
\omega_P$. Since $\omega_P$ is unramified, $c(\pi_P) =
c(\mu)c(\nu) = c(\mu^2)$ will be divisible by $P^2$, and so this
cannot happen. Now part (i) holds because the only special
representation having conductor $P$ (cf. {\it loc. cit.}) is an
unramified twist of the Steinberg representation.

When $\pi$ is dihedral, i.e., when it is automorphically induced
by a character $\Psi$ of a quadratic extension $M/K$, its base
change to $M$ becomes Eisensteinian and so no local component
$\pi_P$ can be special. So (i) implies (ii).

\qed

\bigskip

Recall that the $\ell$-adic representation $V_\ell(\pi_f)$ is free
of rank $16$ over $E \otimes \Q_\ell$, where $E$ is the field of
coefficients of $\pi_f$. Fix an embedding of $\overline \Q$ in
$\overline \Q_\ell$. Then we have
$$
V_\ell(\pi_f) \otimes_{\Q_\ell} \overline \Q_\ell \, \simeq \,
\oplus_{\sigma \in {\rm Hom}(E, \overline \Q)} \, \overline V
_\ell(\pi_f^\sigma),
$$
where each $\overline V_\ell(\pi_f^\sigma)$ is a $16$-dimensional
$\overline \Q_\ell$-representation of ${\mathcal G}_\Q$.

We need to show that for any cuspidal $\pi \in $Coh$_G(C)$, every
Tate class in $\overline V_\ell(\pi_f)$, i.e, a class in
$\overline V_\ell(\pi_f)(2)$ fixed by ${\mathcal G}_\Q$,
contributes to a pole of the $L$-function of $\overline V
_\ell(\pi_f)$. The first object is to analyze the structure of
$\overline V_\ell(\pi_f)$.

\medskip

\noindent{\bf Proposition 8.4} \, \it Let $\tilde K$ denote the
Galois closure of $K$ in $\overline \Q$. There exists a
$2$-dimensional $\overline \Q_\ell$-representation $W_\ell(\pi_f)$
of ${\mathcal G}_K$ such that as ${\mathcal G}_{\tilde K}$-modules
$$
\overline V_\ell(\pi_f) \, \simeq \, \otimes_{\tau \in {\rm
Hom}(K, \overline \Q)} \, W_\ell(\pi_f)^{[\tau]},
$$
where $X^{[\tau]}$ denotes, for any representation $X$ of
${\mathcal G}_{\tilde K}$, the $\tau$-twisted representation
$\alpha \to \tau \circ X(\alpha) \circ \tau^{-1}$.
\rm

\medskip

{\it Proof} \, By a theorem of R.L.~Taylor ([Ta1,2]), proved
independently by D,~Blasius and J.~Rogawski ([B$\ell$-Ro], there
is a $2$-dimensional representation $W_\ell(\pi_f)$ of $\mathcal
G_K$ such that for $S$ a finite set of places of $K$ containing
the archimedean places, the places where $\pi_f$ or $K$ is
ramified, and the ones which divide $\ell$,
$$
L^S(s, \pi) \, = \, L^S(s+1/2, W_\ell(\pi_f)).\leqno(8.5)
$$

Since $K$ is a quartic extension of $\Q$ containing a quadratic
subfield, $\tilde K$ is either $K$ or is a quadratic extension of
$K$. So it makes sense to speak about the base change $\pi_{\tilde
K}$ of $\pi$ to GL$(2)/\tilde K$. It follows that for every
embedding $\tau$ of $K$ in $\overline \Q$, and any extension
$\tilde \tau$ of $\tau$ as an automorphism of $\tilde K$, we have
$$
L^{\tilde S}(s, \pi_{\tilde K} \circ \tau) \, = \, L^{\tilde
S}(s+1/2, W_\ell(\pi_f)^{[\tau]}).
$$
Here $\tilde S$ denotes the set of places of $\tilde K$ above $S$
along with the places where $\tilde K/K$ is ramified.

On the other hand, by definition of the Asai $L$-function of $\pi$
(see section 6), one gets immediately that
$$
L^S(s, \pi_{\tilde K}; r_{K/\Q}) \, = \, L^S(s, \prod_{\tau \in
{\rm Hom}(K, \overline \Q)} \, \pi_{\tilde K} \circ \tau),
\leqno(8.6)
$$
where the formally defined Euler product on the right identifies,
in view of (8.5), with the $L$-function
$$
L^S(s+2, \otimes_{\tau \in {\rm Hom}(K, \overline \Q)} \,
W_\ell(\pi_f)^{[\tau]}).
$$

In view of the Tchebotarev density theorem, it now suffices to
know that
$$
L^S(s+2, \overline V_\ell(\pi_f)) \, = \, L^S(s, \pi_{\tilde K};
r_{K/\Q}). \leqno(8.7)
$$
In fact it is enough to know this over $\tilde K$. In any case,
(8.7) is known by a theorem of Brylinski and Labesse ([BrL]) if we
suitably expand $S$ to a bigger finite set. This suffices for our
purposes. However, it should be pointed out that the need to
expand $S$ is unnecessary thanks to some recent work ([B$\ell$])
of D.~Blasius; to be precise, he has proved the identity (8.7)
with $S$ being just the primes dividing the conductor of $\pi$. We
are now done with the proof of Proposition 8.3.

\qed

\bigskip

Fix an embedding $\sigma$ of $E$ in $\overline \Q_\ell$, and write
$\overline V_\ell(\pi_f)$ for $\overline V_\ell(\pi_f^\sigma)$.
For any number field $k$, let $r_{\ell,k}(\pi_f)$, resp. $r_{{\rm
an}, k}(\pi_f)$, denote the dimension of the $k$-rational Tate
classes in $\overline V_\ell(\pi_f)$, resp. the order of pole at
the $L(s, \overline V_\ell(\pi_f)/k)$ at the edge $s=3$.

\medskip

\noindent{\bf Proposition 8.8} \, \it If $\pi$ is of CM type, we
have for any $k$,
$$
r_{\ell, k}(\pi_f) \, = \, r_{{\rm an}, k}(\pi_f).
$$
\rm

\medskip

{\it Proof}. \, If $\pi$ is of CM type, then $\overline
V_\ell(\pi_f)$ is given by a representation of the Weil group
$W_k$, and the identity follows from the results of [De2]. Done.

\medskip

\noindent{\bf Proposition 8.9} \, When $\pi$ is non-CM, all the
Tate classes in $\overline V_\ell(\pi_f)$ are rational over an
abelian number field $k$, with
$$
r_{\ell, k}(\pi_f) \, \leq \, 2.
$$
\rm

\medskip

{\it Proof}. \, Let $M \supset K$ be a number field. As $\pi$ is
non-CM, the $2$-dimensional Galois representation $W_\ell(\pi_f)$
remains irreducible upon restriction to ${\mathcal G}_M$. Now we
make the following

\medskip

\noindent{\bf Lemma 8.10} \, \it Let $\pi$ be non-CM. Then there
can be a Tate class over $M$, i.e., the character $\chi_\ell^2$
can appear in $\overline V_\ell(\pi_f)$, viewed as a ${\mathcal
G}_M$-module, iff there is a partition
$$
{\rm Hom}(K, \overline \Q) \, = \, \{1, \theta\} \cup \{\tau,
\eta\}, \leqno(\ast)
$$
and an $\ell$-adic character $\lambda_\ell$ of ${\mathcal G}_M$
such that
$$
\lambda_\ell \, \subset \, W_\ell(\pi_f) \otimes
W_\ell(\pi_f)^{[\theta]}, $$
and
$$
\chi_\ell^2\lambda_\ell^{-1} \, \subset \, W_\ell(\pi_f)^{[\tau]}
\otimes W_\ell(\pi_f)^{[\eta]}.
$$
\rm

\medskip

{\it Proof of Lemma 8.10}. \, By Proposition 8.4, we see that a
Tate class exists over $M$ iff
$$
\chi_\ell^2 \, \subset \, \otimes_{\tau \in {\rm Hom}(K, \overline
\Q)} \, W_\ell(\pi_f)^{[\tau]}.
$$
The {\it if} part of the Lemma is now clear, and let us assume
that there is a Tate class over $M$ (to prove the {it only if}
part). We first claim that if $\tau, \beta \in {\rm Hom}(K,
\overline \Q)$ with $\alpha \ne \beta$, the ${\mathcal
G}_M$-module
$$
Y_\ell(\alpha, \beta) \, \simeq \, W_\ell(\pi_f)^{[\alpha]}
\otimes W_\ell(\pi_f)^{[\beta]}
$$
is reducible iff
$W_\ell(\pi_f)^{[\alpha]}$ is a twist, by an $\ell$-adic
character, of $W_\ell(\pi_f)^{[\beta]}$. Indeed, if this fails, we
must have a decomposition (over $M$) of the form
$$
Y_\ell(\alpha, \beta) \, \simeq \, Z_\ell \oplus Z_\ell',
$$
with $Z_\ell, Z_\ell'$ irreducible of dimension $2$. Then the
exterior square of $Y_\ell(\alpha, \beta)$ contains the line
$\Lambda^2(Z_\ell)$. But we also have
$$
\Lambda^2(Y_\ell(\alpha, \beta)) \, \simeq \, ({\rm
sym}^2(W_\ell(\pi_f)^{[\alpha]})\otimes \omega_\ell^\beta) \otimes
({\rm sym}^2(W_\ell(\pi_f)^{[\beta]})\otimes \omega_\ell^\alpha),
$$
where $\omega_\ell^\alpha = \omega^\alpha\chi_\ell$ denotes the
determinant of $W_\ell(\pi_f)^{[\alpha]}$. But as $\pi$ is non-CM,
${\rm sym}^2(W_\ell(\pi_f)^{[\alpha]})$ is irreducible upon
restriction to any open subgroup (such as ${\mathcal G}_M$) of
${\mathcal G}_{K^\alpha}$. The leads to the desired contradiction,
and the claim is proved. Now write
$$
{\rm Hom}(K, \overline \Q) \, = \, \{1, \theta, \tau, \eta\}.
$$
Suppose $Y_\ell(1, \theta)$ is irreducible. Then the existence of
the Tate class implies an isomorphism
$$
Y_\ell(1, \theta) \, \simeq Y_\ell(\tau, \eta)^\vee \otimes
\chi_\ell^2,
$$
and we must have, up to interchanging $\{\tau, \eta\}$,
isomorphisms
$$
W_\ell(\pi_f) \simeq {W_\ell(\pi_f)^{[\tau]}}^\vee \otimes
\nu_\ell \quad {\rm and} \quad W_\ell(\pi_f)^{[\theta]} \simeq
{W_\ell(\pi_f)^{[\eta]}}^\vee \otimes \nu_\ell^{-1}\chi_\ell^2,
$$
for a character $\nu_\ell$. This gives the Lemma in this case,
again up to renaming the embeddings of $K$. So we may assume that
$$
W_\ell(\pi_f)^{[\theta]} \, \simeq \, W_\ell(\pi_f) \otimes
\mu_\ell,
$$
for a character $\mu_\ell$, so that
$$
Y_\ell(1, \theta) \, \simeq \, {\rm sym}^2(W_\ell(\pi_f)) \otimes
\mu_\ell \oplus \omega_\ell\mu_\ell.
$$
The existence of the Tate class implies that
$$
{\rm Hom}_{{\mathcal G}_M}(Y_\ell(1,\theta),
Y_\ell(\tau,\eta)^\vee \otimes \chi^2) \, \ne \, 0,
$$
resulting in an isomorphism $W_\ell(\pi_f)^{[\eta]} \simeq
W_\ell(\pi_f)^{[\tau]}\otimes \nu_\ell$, for a suitable character
$\nu_\ell$. So we have
$$
Y_\ell(\tau,\eta)^\vee \otimes \chi^2) \, \simeq \, ({\rm
sym}^2(W_\ell(\pi_f)^{(\tau]})^\vee \otimes
\nu_\ell^{-1}\chi_\ell^2 \, \oplus \,
(\omega_\ell^\tau\nu_\ell)^{-1}\chi_\ell^2.
$$
Then one of the following must happen:
\begin{enumerate}
\item[(i)] \, $\omega_\ell\nu_\ell \, = \,
(\omega_\ell^\tau\nu_\ell)^{-1}\chi_\ell^2$
\item[(ii)] \, ${\rm sym}^2(W_\ell(\pi_f)^{[\tau]}) \, \simeq \,
{\rm sym}^2(W_\ell(\pi_f))\otimes
(\omega_\ell^2\mu_\ell\nu_\ell)^{-1}\chi_\ell^2$
\end{enumerate}
There is nothing to prove when (i) holds, so we may assume the
identity (ii). So over a finite extension $L$ of $M$, ${\rm
sym}^2(W_\ell(\pi_f)^{[\tau]})$ and ${\rm sym}^2(W_\ell(\pi_f))$
are isomorphic. Then, as is well known (see [Ra2] for example,
though the situation is much simpler here),
$W_\ell(\pi_f)^{[\tau]}$ and $W_\ell(\pi_f)$ will be forced to be
twists of each other by a character of ${\mathcal G}_L$. Since the
determinants of these two modules differ by a finite order
character, they become isomorphic over a finite extension $L_1$ of
$M$. It then follows that
$$
W_\ell(\pi_f)^{[\tau]} \, \simeq \, W_\ell(\pi_f) \otimes
\xi_\ell,
$$
for a character $\xi_\ell$ of ${\mathcal G}_M$. So ${\rm
sym}^2(W_\ell(\pi_f)^{[\tau]})$ is isomorphic to ${\rm
sym}^2(W_\ell(\pi_f))\otimes \xi_\ell^2$. Comparing this with (ii)
and remembering that the symmetric squares do not admit any
self-twist, we obtain
$$
\xi_\ell^2\omega_\ell^2\nu_\ell\mu_\ell\chi_\ell^{-2} \, = \, 1.
\leqno(ii')
$$
A different way to encode the existence of the Tate class (over
$M$) is to note that
$$
{\rm Hom}_{{\mathcal G}_M}(Y_\ell(1,\eta), Y_\ell(\theta,
\tau)^\vee \otimes \chi^2) \, \ne \, 0.
$$
Since $W_\ell(\pi_f)^{[\eta]} \simeq W_\ell(\pi_f)^{[\tau]}\otimes
\nu_\ell \simeq W_\ell(\pi_f)\otimes \xi_\ell\nu_\ell$, we get
$$
Y_\ell(1, \eta) \, \simeq \, ({\rm sym}^2(W_\ell(\pi_f)) \otimes
\xi_\ell\nu_\ell) \oplus \omega_\ell\xi_\ell\nu_\ell.
$$
And since $W_\ell(\pi_f)^{[\theta]} \simeq W_\ell(\pi_f)\otimes
\mu_\ell$,
$$
Y_\ell(\theta, \tau)^\vee \otimes \chi_\ell^2 \, \simeq \, ({\rm
sym}^2(W_\ell(\pi_f)) \otimes
(\mu_\ell\xi_\ell\omega_\ell)^{-1}\chi_\ell^2) \oplus
(\mu_\ell\xi_\ell\omega_\ell)^{-1}\chi_\ell^2.
$$
In view of (ii'), the characters appearing in $Y_\ell(1, \eta)$
and in $Y_\ell(\theta, \tau)^\vee \otimes \chi_\ell^2$ are the
same. This proves the Lemma relative to the partition $\{1,
\eta\}\cup\{\theta, \tau\}$ and $\lambda_\ell =
\omega_\ell\xi_\ell\nu_\ell$. The Lemma follows.

\qed

\medskip

{\it Proof of Proposition 8.9} (contd.) \, As $W_\ell(\pi_f)$,
resp. $W_\ell(\pi_f)^{[\theta]}$, is the restriction of a
representation of ${\mathcal G}_K$, resp. ${\mathcal
G}_{K^\theta}$, $\lambda_\ell$ extends to a character of
${\mathcal G}_{\tilde K}$, where $\tilde K$ is the Galois closure
of $K$. Hence every Tate class over $M$ is already defined over an
abelian extension of $\tilde K$.

Since $W_\ell(\pi_f)$ is Hodge-Tate by [B$\ell$-Ro], so is its
conjugate $W_\ell(\pi_f)^{[\theta]}$. Consequently, $\lambda_\ell$
is also Hodge-Tate, locally algebraic, and is therefore attached
to an algebraic Hecke character. Comparing weights, we see that
$$
\lambda_\ell \, = \, \chi_\ell\beta
\leqno(8.11)
$$
for a {\it finite order character} $\beta$.

Since the dual of $W_\ell(\pi_f)$ is its twist by
$(\omega\chi_\ell)^{-1}$, the identities (8.10) imply the
following:
$$
W_\ell(\pi_f)^{[\theta]} \, \simeq \, W_\ell(\pi_f) \otimes
\beta\omega^{-1},
\leqno(8.12)
$$
and
$$
W_\ell(\pi_f)^{[\eta]} \, \simeq \, W_\ell(\pi_f)^{[\tau]} \otimes
\beta^{-1}\omega^{-\tau}.
$$
Consequently, assuming we have a Tate class over $M$, we may write
$$
\overline V_\ell(\pi_f) \, \simeq \, \left({\rm
sym}^2(W_\ell(\pi_f))\otimes \beta\omega^{-1} \oplus
\beta\chi_\ell\right) \otimes \left({\rm
sym}^2(W_\ell(\pi_f)^{\tau})\otimes \beta^{-1}\omega^{-\tau}
\oplus \beta^{-1}\chi_\ell\right).
\leqno(8.13)
$$
Since $\pi$ is not of CM type, the symmetric square of
$W_\ell(\pi_f))$ is irreducible. So there can be a {\it second
Tate class} over $M$ (which is not a multiple of the first one)
iff we have a non-zero ${\mathcal G}_M$-homomorphism
$$
\phi: \chi_\ell^2 \, \rightarrow \, {\rm sym}^2(W_\ell(\pi_f))
\otimes {\rm sym}^2(W_\ell(\pi_f)^{[\tau]}).
\leqno(8.14)
$$
The irreducibility of ${\rm sym}^2(W_\ell(\pi_f))$ implies that
there can be at most one such $\phi$. Consequently, using the fact
any Tate class over a number field $k$ remains a Tate class over
$M = kK$, we get
$$
r_{\ell, k}(\pi_f) \, \leq \, 2. \leqno(8.15)
$$

\medskip

Now we show that all the Tate classes are rational over an {\it
abelian} number field (when $\pi$ is non-CM). If $r_{\ell,
\overline \Q}(\pi_f) = 1$, then the fact that $\overline
V_\ell(\pi_f)$ is a representation of ${\mathcal G}_\Q$ implies
that for a finite order character $\nu$ of ${\mathcal G}_\Q$ which
becomes trivial when restricted to ${\mathcal G}_M$,
$\nu\chi_\ell$ must be a summand of $\overline V_\ell(\pi_f)$.
Consequently, the Tate class is defined over the cyclic extension
of $\Q$ cut out by $\nu$.

It remains to consider the case $r_{\ell, \overline \Q}(\pi_f) =
2$. Here, by the discussion in (a), we must have (8.10) through
(8.14). It follows from (8.14) and the irreducibility of ${\rm
sym}^2(W_\ell(\pi_f))$ that we must have, as ${\mathcal G}_{\tilde
K}$-modules:
$$
{\rm sym}^2(W_\ell(\pi_f)^{[\tau]})) \, \simeq \, {\rm
sym}^2(W_\ell(\pi_f)) \otimes \nu, \leqno(8.16)
$$
Then over a finite extension $L$ (where $\nu$ becomes trivial),
${\rm sym}^2(W_\ell(\pi_f)^{[\tau]}))$ and ${\rm
sym}^2(W_\ell(\pi_f))$ are isomorphic. Then, as seen in the proof
of Lemma 8.10, we must have (as ${\mathcal G}_L$-modules,
$$
W_\ell(\pi_f)^{[\tau]}) \, \simeq \, W_\ell(\pi_f) \otimes \mu,
\leqno(8.17)
$$
for a character $\mu$ of ${\mathcal G}_{\tilde L}$. Since $\pi$ is
non-CM, there can be no character other than $\mu$  occurring in
$W_\ell(\pi_f)^\vee \otimes W_\ell(\pi_f)^{[\tau]})$. Hence $\mu$
extends to a character of ${\mathcal G}_{\tilde K}$, and (8.17) is
valid over $\tilde K$. Now putting Lemma 8.10 and (8.17) together,
we see that conjugation of $W_\ell(\pi_f)$ by the nontrivial
automorphism of $K/F$, call it $\theta$, is equivalent to a twist
of $W_\ell(\pi_f)$. Consequently,
$$
W_\ell(\pi_f) \otimes W_\ell(\pi_f)^{[\theta]} \, \simeq \, ({\rm
sym}^2(W_\ell(\pi_f)) \otimes \alpha) \oplus \alpha\chi_\ell,
\leqno(8.18)
$$
for a finite order character $\alpha$. It follows that
$$
As_{K/F}(W_\ell(\pi_f)) \, \simeq \, \beta_\ell \oplus
\varphi\chi_\ell, \leqno(8.19)
$$
where $\beta_\ell$ is an irreducible $3$-dimensional
representation of ${\mathcal G}_F$ and $\varphi$ is a finite order
character. Then $\beta_\ell$ must be essentially self-dual,
meaning that its symmetric square contains a $1$-dimensional
summand. We claim that there can be no other one-dimensional
summand. Indeed, the restriction of $sym^2(\beta_\ell)$ to
${\mathcal G}_K$ is, up to a twist, isomorphic to
$sym^2(sym^2(W_\ell(\pi_f)))$, which splits as
sym$^4(W_\ell(\pi_f))$ and a character. Now since $W_\ell(\pi_f)$
is irreducible upon restriction to any open subgroup of ${\mathcal
G}_K$, it follows that $sym^4(W_\ell(\pi_f))$ is irreducible,
proving the claim. Now if $\tau$ denotes the non-trivial
automorphism of $F$, we get from Lemma 8.10 and (8.14) the
consequence that $sym^2(W_\ell(\pi_f)^{[\tau]})$ is a twist of
$sym^2(W_\ell(\pi_f))$. Putting these together we get the
decomposition
$$
As_{K/F}(W_\ell(\pi_f)) \otimes As_{K/F}(W_\ell(\pi_f)^{[\tau]})
\, \simeq \, \sigma_\ell \oplus \varphi\varphi^\tau\chi_\ell^2
\oplus \xi\chi_\ell^2, \leqno(8.20)
$$
where $\xi$ is a finite order character, and $\sigma_\ell$ is a
$14$-dimensional representation whose irreducible summands are
$5$-dimensional or $3$-dimensional. Recall that the restriction of
$As_{K/Q}(W_\ell(\pi_f))$ to ${\mathcal G}_F$ is isomorphic to
$As_{K/F}(W_\ell(\pi_f)) \otimes
As_{K/F}(W_\ell(\pi_f)^{[\tau]})$. Clearly Gal$(F/\Q)$ permutes
the irreducible constituents on the right of (8.20), and for
dimension reasons it must preserve the set $\{\varphi\varphi^\tau,
\xi\}$. Since $\varphi\varphi^\tau$ is invariant, it extends to a
character $\varphi'$, say, of ${\mathcal G}_\Q$. This means that
$\xi$ must also be $\tau$-invariant and must extend to a character
$\varphi''$, say, of ${\mathcal G}_\Q$. It follows that the Tate
classes are all defined over the compositum of the cyclic
extensions of $\Q$ cut out by $\varphi'$ and $\varphi''$. Done.

\qed

\bigskip

We can write
$$
r_{\ell, \Q^{\rm ab}}(\pi_f) \, = \, \sum_\nu r_\ell(\pi_f; \nu),
\leqno(8.21)
$$
where $\nu$ runs over the finite order characters of Gal$(\Q^{\rm
ab}/\Q)$ and $r_\ell(\pi_f; \nu)$ denotes the rank of the
$\nu$-isotypic subspace $Ta_\ell(\pi_f; \nu)$ of the space of
$\pi_f$-Tate classes over $\Q^{\rm ab}$. Of course $r_\ell(\pi_f;
1) \, = \, r_{\ell, \Q}(\pi_f)$. Thanks to Proposition 8.9, we
know (for $\pi$ non-CM) that
$$
r_\ell(\pi_f; \nu) \, \leq \, r_{\ell, \Q^{\rm ab}}(\pi_f) \, \leq
2.
$$

\medskip

\noindent{\bf Proposition 8.22} \, \it Let $K/\Q$ be Galois, and
$\pi$ non-CM.
\begin{enumerate}
\item[(a)] \, $r_{\ell, \Q^{\rm ab}}(\pi_f)$ is non-zero iff a
twist of $\pi$ is a base change from a quadratic subextension of
$K$.
\item[(b)] \, $r_{\ell, \Q^{\rm ab}}(\pi_f)$ equals $2$ iff a
twist of $\pi$ is a base change from $\Q$.
\item[(c)] \, The
following are equivalent:
\begin{enumerate}
\item[(i)] \, $r_\ell(\pi_f; \nu) \, = \, 2$ for some $\nu$.
\item[(ii)] \, A twist of $\pi$ is a base change from $\Q$, and
$K/\Q$ is biquadratic.
\end{enumerate}
\end{enumerate}
\rm

\medskip

{\it Proof}. \, (a) \, Suppose $r_{\ell}(\pi_f, \nu) \ne 0$ for a
character $\nu$ of the absolute Galois group of $\Q$. Let $\tilde
\nu$ be a character of the idele class group of $K$ which extends
the idele class character of $\Q$ attached to $\nu$ by class field
theory. Then we know that
$$
\overline V(\pi_f \otimes \tilde \nu^{-1}) \, \simeq \, \overline
V(\pi_f) \otimes \nu^{-1},
$$
so that
$$
r_\ell(\pi_f,\nu) \ne 0 \, \Leftrightarrow \, r_\ell(\pi_f \otimes
\tilde\nu) \ne 0. \leqno(8.23)
$$
So we may assume, by replacing $\pi$ be $\pi \otimes
\tilde\nu^{-1}$, that $r_{\ell, \Q}(\pi_f)$ is itself non-zero,
i.e., that $\chi_\ell^2$ appears in the $\mathcal G_\Q$-module
$\overline V_\ell(\pi_f)$. Then by Lemma 8.10, there is a
partition $\{1,\theta\}\cup\{\tau,\eta\}$ of Hom$(K,\overline
\Q)$, and an $\ell$-adic character $\lambda_\ell$ of the form
$\beta\chi_\ell$, with $\beta$ of finite order, such that $(\ast)$
holds, with $\beta\beta^\tau = 1$. Then in the unitary
normalization,
$$
\beta \, \subset \, \pi \boxtimes \pi^\theta. \leqno(8.24)
$$
We have to show that a twist of $\pi$ is a base change from a
quadratic subfield.

First consider the {\it biquadratic case}, when $\theta^2 = \tau^2
= 1$ and $\eta = \theta\tau$. Then $\pi \boxtimes \pi^\theta$ is
$\theta$-invariant, and so $\beta^\theta$ is contained in $\pi
\boxtimes \pi^\theta$. Put $\mu = \beta^{-1}\beta^\theta$. Then
$$
\pi^\theta \simeq \pi^\vee \otimes \beta \simeq \pi^\vee \otimes
\beta^\theta, \quad {\rm so} \quad \pi^\vee \otimes \mu \simeq
\pi^\vee. \leqno(8.25)
$$
Consequently $\pi$ admits a self-twist by $\mu^{-1}$, which is
impossible as we are in the non-CM case, unless $\mu =1$. Then
$\beta = \beta^\theta$. Hence $\beta$ comes from a character
$\gamma$, say, of the quadratic subfield $F$ fixed by $\theta$.
Then $As_{K/F}(\pi)$ contains $\gamma$ as an isobaric summand. If
$\tilde \gamma$ is an idele class character of $K$ restricting to
$\gamma$, we know that $As_{K/F}(\pi \otimes \tilde \gamma^{-1})
\simeq As_{K/F})(\pi) \otimes \nu^{-1}$, which contains $1$. So a
twist of $\pi$ is a base change from $F$.

Now let $K/\Q$ be {\it cyclic}. If $\theta$ has order $2$, the
argument of the biquadratic case carries through. So let $\theta$
have order $4$, so that $\tau = \theta^2$ and $\eta = \tau\theta =
\theta^3$. Now $\beta$ occurring in $\pi \boxtimes \pi^\theta$
implies that $\pi^\tau$ is $\pi^\vee \otimes \mu$, with $\mu =
\beta^\theta/\beta$. In other words, $\mu$ occurs in $\pi\boxtimes
\pi^\tau$, and the rest of the argument is similar.

Conversely, suppose a twist of $\pi$ is a base change from a
quadratic subfield $F$ of $K$. Then $As_{K/F}(\pi)$ contains a
character $\xi$ as an isobaric summand. Then the transfer to $\Q$
of the Galois character defined by $\xi$ is necessarily of the
form $\chi_\ell^2\nu$ and occurs in the $16$-dimensional
Gal$(\overline \Q/\Q)$-representation, which corresponds to
$As_{F/\Q}(As_{K/F}(\pi))$. In other words $r_\ell(\pi_f,\nu) \ne
0$. We are done with this part.

\medskip

(b) \, Let $r_\ell(\pi_f, \nu) = 2$ for some $\nu$. Let $F$ be the
quadratic extension of $\Q$ contained in $K$, given by part (a),
such that a twist of $\pi$ is the base change of a cusp form
$\pi_0$ on GL$(2)/F$ of central character $\omega_0$. Note that
the restriction of $\alpha$ to $F$ will necessarily be the
non-trivial automorphism of $F$. We may replace $\pi$ by the
appropriate twist and assume that it is just $\pi_{0,K}$. It is
easy to see that
$$
As_{K/F}(\pi_{0,K}) \, \simeq \, sym^2(\pi_0) \boxplus
\omega_0\delta, \leqno(8.26)
$$
where $\delta = \delta_{K/F}$. Since a second character occurs in
the Galois module, there must exist a character $\varphi$ of $F$
such that
$$
sym^2(\pi_0^\alpha) \, \simeq \, sym^2(\pi_0^\vee) \otimes \varphi
\, \simeq sym^2(\pi_0) \otimes \varphi\omega_0^{-2}. \leqno(8.27)
$$
If we put $\xi = \varphi^\alpha/\varphi$, we see that
$$
sym^2(\pi_0) \, \simeq \, sym^2(\pi_0) \otimes \xi.
$$
Comparing central characters, $\xi^3 = 1$. If $\xi$ is
non-trivial, $\pi_0$ will become dihedral over a cubic extension
of $F$ without being so already over $F$, and this is not
possible. So $\xi = 1$, and $\varphi$ extends to a character of
$\Q$. Furthermore, since the corresponding Tate class over $\Q$
forces the restriction of $\varphi\omega_0^{-2}$ to $\Q$ to be
trivial. So we may write $\varphi\omega_0^{-2}$ as
$\gamma/\gamma^\alpha$, for some character $\gamma$ of $F$. Then
$sym^2(\pi_0) \otimes \gamma$ is a base change from $\Q$.

We {\it claim} that $\pi_0^\alpha$ and $\pi_0$ are themselves
twist equivalent. Indeed, if $L = F(\gamma)$, the symmetric
squares of $\pi_{0,L}^\alpha$ and $\pi_{0,L}$ are equivalent.
Thus, as seen above in the proof of Lemma 8.10, we must have
$\pi_{0,L}^\alpha \simeq \pi_{0,L} \otimes \psi$, for a character
$\psi$ of $L$. Then $\pi_0^\alpha$ and $\pi_0$ become isomorphic
after base change to the extension $R := L(\psi)$, which implies
the claim.

So some character $\beta$ of $F$ appears in $\pi_0 \boxtimes
\pi_0^\alpha$. Now the rest of the proof goes as in the proof of
part (a) following (8.24).

It remains to prove the converse.

\medskip

\noindent{\bf Proposition 8.28} \, \it Let $\pi'$ be a cusp form
on GL$(2)/\Q$ of weight $2$ and central character $\omega'$, and
let $K$ be a quartic totally real extension containing a quadratic
subfield $F$. Then $As_{K/\Q}(\pi'_K)$ is automorphic with the
decomposition
$$
As_{K/\Q}(\pi'_K) \, \simeq \, sym^4(\pi') \boxplus (sym^2(\pi')
\otimes I_F^\Q(\delta) \otimes \omega'\epsilon) \boxplus
(sym^2(\pi') \otimes \omega'\epsilon) \boxplus {\omega'}^2
\boxplus {\omega'}^2\delta_1,
$$
where $\delta$ is the quadratic character of $F$ associated to $K$
with restriction $\delta_1$ to $\Q$ and $\epsilon$ is the
quadratic character of $\Q$ attached to $F/\Q$. \rm

\medskip

It is easy to see that this Proposition concludes the proof of
part (b) of Proposition 8.2.2.

\medskip

{\it Proof of Proposition 8.28}. \, In what follows we will do
everything formally and treat all the representations as
admissible representations. Once the identity is proved, however,
the fact that the right hand side is automorphic (thanks to Kim
[K]) will imply the automorphy of the left hand side. To be
precise, we do not really need the automorphy, only the ability to
understand the behavior of the relevant $L$-function at the right
edge.

An immediate consequence of the Asai representation is the
following identity:
$$
As_{K/F}(\pi'_K) \boxplus I_K^F(sym^2(\pi'_K)) \, \simeq \,
sym^2(I_K^F(\pi'_K)). \leqno(8.29)
$$
Since the Asai construction is compatible with doing it in stages,
the obvious analogue of (8.29) furnishes the isomorphism
$$
As_{K/\Q}(\pi'_K) \boxplus I_F^\Q(sym^2(sym^2(\pi'_F) \boxplus
\omega'_F\delta)) \, \simeq \, sym^2(I_{F}^\Q(sym^2(\pi'_F))
\boxplus I_{F}^\Q(\omega'_F\delta)). \leqno(8.30)
$$
Since $I_F^\Q(sym^2(\pi'_F))$ (resp. $I_F^\Q(\omega'_F\delta)$)
identifies with $sym^2(\pi')\boxtimes (1\boxplus \epsilon)$ (resp.
$\omega' I_F^\Q(\delta)$, and since $sym^2(I_F^\Q(\delta))$ is
just $1 \boxplus \epsilon \boxplus \delta_1$, the right hand side
of (8.30) becomes
$$
sym^2(sym^2(\pi') \boxtimes(1 \boxplus \epsilon)) \boxplus
({\omega'}^2\otimes(1\boxplus \epsilon \boxplus \delta_1))
\boxplus (sym^2(\pi') \otimes(1 \boxplus \epsilon) \boxtimes\omega
I_F^\Q(\delta)).
$$
Moreover, $I_F^\Q(sym^2(sym^2(\pi'_F))$ identifies with
$$
sym^2(sym^2(\pi')) \boxplus (sym^2(sym^2(\pi')) \otimes \epsilon)
\boxplus {\omega'}^2\boxtimes(1\boxplus \epsilon) \boxplus
(sym^2(\pi') \otimes \omega'\otimes I_F^\Q(\delta)).
$$
But
$$
sym^2(sym^2(\pi') \boxtimes(1 \boxplus \epsilon) \, \simeq \,
(sym^2(sym^2(\pi')) \boxtimes sym^2(1\boxplus\epsilon)) \boxplus
(\Lambda^2(sym^2(\pi')) \boxtimes \Lambda^2(1\boxplus\epsilon)),
$$
which simplifies as
$$
sym^2(sym^2(\pi') \boxtimes(1 \boxplus \epsilon \boxplus 1)
\boxplus (sym^2(\pi') \otimes \omega') \otimes \epsilon.
$$
Consequently,
$$
As_{K/\Q}(\pi'_K) \, \simeq \, sym^2(sym^2(\pi')) \boxplus
(sym^2(\pi') \otimes \omega'\epsilon) \boxplus {\omega'}^2\delta_0
\boxplus (sym^2(\pi') \boxtimes \epsilon\omega' \boxtimes
I_F^\Q(\delta)). \leqno(8.31)
$$
The Proposition now follows thanks to the identity
$$
sym^2(sym^2(\pi')) \, \simeq \, sym^4(\pi') \boxplus {\omega'}^2.
$$
Done.

\qed

\bigskip

{\it Proof of Proposition 8.22} (contd.) \,

\noindent{(c)} \, In view of part (b), it suffices to show that if
$\pi$ is a twist of the base change $(\pi_1)_K$ of a weight $2$
cusp form $\pi_1$ on GL$(2)/\Q$, $r_\ell(\pi_f;\nu) =2$ for some
$\nu$ iff $K/\Q$ is biquadratic. Thanks to Proposition 8.28, it is
then enough to show that
$$
\delta_1 \, = \, 1 \, \Leftrightarrow \, K/\Q \quad {\rm
biquadratic}.
$$
By definition, $\delta_1$ is the restriction to $\Q$ of $\delta =
\delta_{K/F}$. When $K$ is biquadratic, we have $K = LF$ for a
quadratic extension $L/\Q$, and if $\delta'$ denotes the quadratic
character of $\Q$ attached to $L$, $\delta = \delta' \circ
N_{F/\Q}$. This gives what we want in this case as for any $x$ in
$I_\Q/\Q^\ast$,
$$
\delta_1(x) \, = \, \delta'(x^2) \, = \, 1.
$$
So we may suppose that $K/\Q$ is cyclic. Then if $\xi$ is the
quartic character of $\Q$ attached to $K/\Q$, $\delta = \xi\circ
N_{F/\Q}$, so that on the idele classes of $\Q$, $\delta_1 = \xi^2
\ne 1$. Done.

\qed

\bigskip

{\it Proof of Proposition 8.2}. \, Thanks to Proposition 8.8, we
may assume that $\pi$ is not of CM type. The proof of part (a)
(resp. (b)) of Proposition 8.22 shows in fact that $r_{{\rm an},
\Q^{\rm ab}}(\pi_f)$ is $\geq 1$ (resp. $= 2$) iff a twist of
$\pi$ is a base change from a quadratic subfield (resp. from
$\Q$).  The desired equality of $r_{{\rm an}, \Q^{\rm ab}}(\pi_f)$
and $r_{\ell, \Q^{\rm ab}}(\pi_f)$ then follows from (8.26) and
Proposition 8.28, modulo the well known fact that for any isobaric
automorphic representation $\Pi$ on any GL$(n)$, $L(s, \Pi)$ is
non-vanishing at the right edge and has no pole there if $1$ does
not occur in the isobaric decomposition of $\Pi$. We are also
using here the modularity results of [Ra3] and [K].

\qed

\medskip

As noted already, Proposition 8.1 follows from Proposition 8.2. We
are now done with this section.

\bigskip

\noindent{\bf Remark 8.33}: \, Suppose $K/\Q$ is {\it non-normal}
with an intermediate field $F$, and $\pi$ a twist of the base
change of a form on GL$(2)/\Q$. Then Proposition 8.28 still holds,
and moreover, it is not hard to see that $\delta_1 \ne 1$. So
$r_\ell(\pi_f; \nu)$ is never $2$ for such a $\pi$. On the other
hand, there are $\pi$ in the non-Galois case with $r_{\rm
an}(\pi_f) \ne 0$, but with $\Pi:= As_{K/F}(\pi)$ cuspidal. Since
$L(s, As_{F/\Q}(\Pi))$ has a pole at $s=1$, its period integral
over $Z(\A_\Q){\rm GL}(4, \Q)\backslash {\rm GL}(4, \A_\Q)$ is
non-zero. This should mean, by a conjecture of Jacquet, that $\Pi$
comes from a unitary group in four variables associated to $F/\Q$.
It is not clear (to this author) as to how to use it to deduce the
Tate conjecture in that case, whence the Galois assumption in the
second half of Theorem A.

\medskip

\vskip 0.2in

\section{Matching Poles with algebraic cycles}

\bigskip

We can define $r_{\rm alg}(\pi_f; \nu)$ in the obvious way, and
the results of the previous section says that
$$
r_{\rm alg}(\pi_f; \nu) \, \leq \, r_{\ell}(\pi_f; \nu) \, = \,
r_{\rm an}(\pi_f; \nu) \, \leq \, 2. \leqno(9.1)
$$

Theorem A$^\prime$ (and hence Theorem A) will follow once we
establish the following

\medskip

\noindent{\bf Proposition 9.2} \, \it Let $K$ be a quartic,
Galois, totally real number field, $C = C_0({\mathfrak N})$, and
$\pi$ a non-CM cusp form on GL$(2)/K$ of weight $2$ contributing
to Coh$_C$. Then for any Dirichlet character $\nu$, we have
$$
r_{\rm alg}(\pi_f; \nu) \, = \, r_{\rm an}(\pi_f; \nu).
$$
\rm

\medskip

{\it Proof}. \, Suppose $r_{\rm an}(\pi_f; \nu) = 1$ for some
$\nu$. Then by part (a) of Proposition 8.22, there is a quadratic
subfield $F$ of $K$, a cusp form $\pi_0$ on GL$(2)/F$ with central
character $\omega_0$, and a character $\xi$ of $K$, such that
$$
\pi \, \simeq (\pi_0)_K \otimes \xi.
$$
Let $\xi_0$ denote the restriction of $\xi$ to (the idele classes
of) $F$. Then we get (using (8.26))
$$
As_{K/F}(\pi) \, \simeq \, As_{K/F}((\pi_0)_K) \otimes \xi_0 \,
\simeq \, (sym^2(\pi_0) \otimes \xi_0) \boxplus
\omega_0\xi_0\delta, \leqno(9.3)
$$
with $\delta = \delta_{K/F}$. If $\mu_1$ denotes the restriction
to $\Q$ of $\mu_0: = \omega_0\xi_0\delta$, then it occurs in
$As_{K/\Q}(\pi)$. In other words, $\nu = \overline\mu_1$ and
$$
\mu_0\chi_\ell^2 \, \subset \, \overline V_\ell(\pi_f)
$$
as Gal$(\overline \Q/\Q)$-modules. Choose a finite order character
$\mu$ of $K$ with restriction $\mu_0$ to $F$. Then $r_{\rm
alg}(\pi_f; \overline\mu_1) \, = \, r_{\rm alg}(\pi_f \otimes
\overline\mu)$, we need only show that
$$
r_{\rm alg}(\pi_f \otimes \overline\mu) \, \ne \, 0. \leqno(9.4)
$$
But (9.3) implies that $L(s, As_{K/F}(\pi \otimes \overline \mu))$
has a pole at the right edge, which implies by the residue formula
of [HLR],
$$
\int\limits_{{\rm GL}(2,F)Z(\A_F)\backslash {\rm GL}(2, \A_F)} \,
\phi(g)\overline\mu({\rm det}(g)) dg \ \ne \, 0, \leqno(9.5)
$$
for some function $\phi$ in the space of $\pi$, with $Z$ denoting
the center of GL$(2)$. In other words, the integral of a
$(2,2)$-form $\eta_\phi$ on the Hilbert modular fourfold defined
by $\phi$ has non-zero $\overline \mu$-twisted period over a Hecke
translate of the embedded Hilbert modular surface attached to $F$.
Look at the corresponding twisting self-correspondence of the
fourfold (see the end of section 5), which defines (for a suitable
$g_f$ in $G_f$) a $\overline\mu$-twisted H-Z cycle
$Z({\overline\mu}) = {}^FZ^\ast_{g_f, C}(\overline\mu)$ (of
codimension $2$). And we obtain
$$
\int_{Z(\overline\mu)} \, \eta_\phi \, \ne \, 0,
$$
proving (9.4).

\medskip

So we may assume that $r_{\rm an}(\pi_f; \nu) = 2$, which implies,
thanks to part (c) of Proposition 8.22, that $K/\Q$ is biquadratic
and that
$$
\pi \, \simeq (\pi')_K \otimes \psi,
$$
for a cusp form $\pi'$ on GL$(2)/\Q$ of central character
$\omega'$, and a character $\psi$ of $K$ (with restriction
$\psi_1$ to $\Q$). Applying Proposition 8.28, we get the embedding
$$
\psi_1{\omega'}^2\chi_\ell^2 \oplus \psi_1{\omega'}^2\chi_\ell^2
\, \subset \, \overline V_\ell(\pi_f). \leqno(9.6)
$$
Let $\mu$ be a character of $K$ with restriction $\mu_1: =
\psi_1{\omega'}^2$ to $\Q$. Then again $\nu = \overline\mu_1$, and
we need to show that
$$
r_{\rm alg}(\pi_f \otimes \overline\mu) \, = \, 2.
\leqno(9.7)
$$
Choose two quadratic subfields $F, E$, say, of $K$. Then $\pi
\otimes \overline\mu$ is a base change from both $F$ and $E$. So
we get two twisted, codimension $2$ algebraic cycles $Z: =
{}^FZ(\overline\mu; g_f)$ and $Z':=
{}^EZ^\ast_{g'_f,C}(\overline\mu)$ on $\tilde S$, for suitable
$g_f, g'_f \in G_f$. These are {\it homologically non-trivial}
because the period integrals of suitable $(2,2)$-forms over these
cycles are non-zero, the reason being that they arise as the
residues of the associated degree four Asai $L$-functions.

But these two cycles may be proportional in the $\pi_f$-component
of the cohomology. So we have to replace one of them with a
suitable twisted version.

\medskip

Let $\theta, \alpha$ be the automorphisms of $K$ with respective
fixed fields $F, E$. Then the restriction of $\alpha$ to $F$ is
the non-trivial automorphism of $F$. Fix an embedding $w: K
\hookrightarrow \R$ and order the archimedean places of $K$ as
$(w, \alpha w, \theta w, \alpha\theta w)$. Given any signature
distribution $s = (s_1, s_2, s_3, s_4)$, with each $s_j$ being $+$
or $-$, we can define a {\it real analytic automorphism} $\tau_s$
of $K \otimes \C - K \otimes \R \simeq \C^4-\R^4$ by
$$
\tau_s(z_1,z_2,z_3,z_4) \, = \, (\tau_{s_1}(z_1), \tau_{s_2}(z_2),
\tau_{s_3}(z_3), \tau_{s_4}(z_4)) \leqno(9.8)
$$
where $\tau_{s_j}$ is the identity, resp. complex conjugation, if
$s_j$ is $+$, resp. $-$. Each such involution acts on the Hilbert
modular fourfold $\tilde S_C$ and its rational Betti cohomology.
It also commutes with the Hecke action and we obtain a
decomposition
$$
V_B(\pi_f) \, \simeq \, \oplus_{s \in \Sigma} V_B(\pi_f)^s,
\leqno(9.9)
$$
where $\Sigma$ runs over all the signature distributions, and
$V_B(\pi_f)^s$ denotes the $s$-eigenspace of $V_B(\pi_f)$. Note
that $\Sigma$ forms a group under componentwise multiplication
with identity $(+,+,+,+)$. Also, each eigenspace $V_B(\pi_f)^s$ is
one-dimensional.

\medskip

If $\xi$ is a finite order character of $K$, its component at any
archimedean place will be $1$ or the {\it sign character}, and we
will define the {\it signature} $s(\xi)$ of $\xi$ to be
$(s(\xi)_1, s(\xi)_2, s(\xi)_3, s(\xi)_4)$, where $s(\xi)_j$ is
the sign of $\xi_w$, resp. $\xi_{\tau w}$, resp. $\xi_{\theta w}$,
resp. $\xi_{\theta\tau w}$, for $j=1$, resp. $j=2$, $j=3$, $j=4$.
(At any archimedean place $u$, the sign of $\xi_u$ is $+$, resp.
$-$, if $\xi_u$ is trivial, resp. non-trivial.) It is easy to
check the following:

\medskip

\noindent{\bf Lemma 9.10} \, \it The twisting correspondence
$R(\xi)$ sends any vector in $V_B(\pi_f)^s$ into $V_B(\pi_f
\otimes \xi)^{ss(\xi)}$. \rm

\medskip

\noindent{\bf Lemma 9.11} \, \it There exists a finite order
character $\xi$ of $K$ such that
\begin{enumerate}
\item[(i)] \, $s(\xi) \, = \, (+,+,-,-)$
\item[(ii)] \, $\xi|_E \, = \, 1$.
\end{enumerate}
\rm

\medskip

{\it Proof of Lemma 9.11}. \, Pick any finite order character
$\lambda$ of $K$ of signature $(+,+,+,-)$. Then $\lambda^\alpha$
has signature $(+,+,-,+)$. Put $\xi = \lambda/\lambda^\alpha$.
Then $s(\xi) = (+,+,-,-)$. And $\xi$ also satisfies (i) because
$E$ is the fixed field of $\alpha$. Done.

\medskip

{\it Proof of Proposition 9.2} (contd.) \, Pick a $\xi$ as in the
Lemma. As $\xi|_\Q = 1$,
$$
r_{\rm an}(\pi_f \otimes \xi; \nu) \, = \, r_{\rm an}(\pi_f; \nu).
\leqno(9.12)
$$
And since $\pi \simeq (\pi')_K \otimes \mu$, we have
$$
As_{K/F}(\pi \otimes \xi\overline\mu) \, \simeq \,
sym^2(\pi'_F)\otimes \xi|_F \boxplus \xi|_F \leqno(9.13)
$$
and (since $\xi|_E = 1$)
$$
As_{K/E}(\pi \otimes \xi\overline\mu) \, \simeq \, sym^2(\pi'_E)
\boxplus 1. \leqno(9.14)
$$
So $L(s, As_{K/E}(\pi\otimes\xi\overline\mu))$ has a simple pole
at the edge. (But $L(s, As_{K/F}(\pi\otimes\xi\overline\mu))$ has
no pole at the right edge as $\xi|_F$ is non-trivial by
construction.) Consequently, if we put
$$
Z'': = \, {}^EZ^\ast_{g''_f, C_1}(\overline\mu\xi),
$$
we have (for suitable $g_f'' \in G_f$ and compact open $C_1$,
$$
\int_{Z''} \eta_\phi \, \ne \, 0, \leqno(9.15)
$$
for some $\phi$ in the space of $\pi$.

\medskip

\noindent{\bf Lemma 9.16} \, \it The space spanned by the classes
of $Z$, $Z'$ and $Z''$ in $V_B(\pi_f)$ has dimension $2$. \rm

\medskip

{\it Proof}. \, As we have seen, these three cycles are all
homologically non-trivial in the $\pi_f$-component. If $[Z]$,
$[Z']$ are not proportional, there is nothing to prove. So we may
suppose that they span a line $L$, say, in $V_B(\pi_f)$. Since $Z$
(resp. $Z'$) comes from $F$ (resp. $E$), $[Z]$ (resp. $[Z']$) has
a non-zero component in $V_B(\pi_f)^s$ for some
$s=(s_1,s_2,s_3,s_4)$ iff $s_1=s_3$ and $s_2=s_3$ (resp. $s_1=s_2$
and $s_3=s_4$). So $L$ lies in $Y:= V_B(\pi_f)^{(+,+,+,+)} \oplus
V_B(\pi_f)^{(-,-,-,-)}$. Now since $\xi$ has signature
$(+,+,-,-)$, $[Z'']$ cannot lie in $Y$, thanks to  Lemma 9.10.
Done.

\qed

\medskip

The Proposition is now proved, as is Theorem A$^\prime$, which
implies Theorem A.

\qed

\medskip

\noindent{\bf Remark 9.17}: \, The referee has suggested the
following clever, alternate approach to proving that the algebraic
cycles in $V_B(\pi_f)$ span a plane: The fact that the period
integral over $Z$ (resp. $Z'$) is non-zero implies that it gives
rise to an SL$(2, \A_F)$-invariant (resp. SL$(2, \A_E)$-invariant)
linear form on the space of $\pi$. So if the homology classes of
$Z$ and $Z'$ are proportional, then we would get an SL$(2,
\A_K)$-invariant linear form on the space of $\pi$, which is
impossible as $\pi$ has no non-trivial intertwining map into an
abelian representation. Hence $[Z], [Z']$ are not proportional.
Since they are non-trivial, they must span a $2$-dimensional
vector subspace of $V_B(\pi_f)$.

\bigskip

We will now give a {\it justification of the remark coming right
after the statement of Theorem A} in the Introduction. We have to
show that there are codimension $2$ cycles which are not
intersections of divisors. We {\it claim} that this is {\it
always} the case for classes coming from cuspidal cohomology,
i.e., that the intersection of divisors will never hit
$V_B(\pi_f)$ for any cuspidal $\pi_f$. Indeed, we can ignore those
divisors supported on the boundary, and for the others, the
representations of $G(\R)$ contributing to $IH^2(S_C^\ast(\C),
\C)$ are one-dimensional, and the tensor product of two such
cannot contain any discrete series representation, which is what
contributes to $H^4_{\rm cusp} \subset IH^4(S_C^\ast(\C), \C)$.
Done.

\qed

\bigskip

\noindent{\bf Remark 9.18} \, It is time to make some comments on
the {\it CM situation} and raise a few questions. Let $\pi$ be a
cohomological cusp form of CM type, say of trivial central
character. Then $\pi$ is the automorphic induction of the unitary
version $\Psi^u$ of a weight one Hecke character $\Psi$ of a CM
quadratic extension $M$ of $K$. Let us take $K$ to be Galois over
$\Q$. If $\tilde M$ denotes the Galois closure of $M$, then the
restriction of $V_\ell(\pi_f)$ to ${\mathcal G}_{\tilde M}$ splits
as a direct sum of sixteen $1$-dimensional representations, six of
which could be of Tate type. For any of these $1$-dimensionals to
define a Tate class (of codimension $2$) over an extension field,
however, the infinity type of the corresponding Hecke character
must be that of the square of norm (or its inverse depending on
the normalization). Looking at the different possible CM types one
sees that this is possible only if $K$ contains a quadratic
subfield, say $F$, such that $M/F$ is biquadratic. Let us analyze
this case when $\pi$ is a base change to $K$ of a CM cusp form
$\pi_0$ of $F$, with $\pi_0$ being defined by a weight one Hecke
character $\varphi$ of a CM quadratic extension $L$ of $F$ so that
$M = LK$. Let $\delta$ (resp. $\nu$) denote the quadratic
character of $F$ attached to $K$ (resp. $L$), and let $\varphi_0$
be the restriction of $\varphi$ to $F$ (which corresponds to the
transfer of the Galois character defined by $\varphi$). Then
$\pi_0$ has central character $\varphi_0\nu$ and its symmetric
square is $I_L^F((\varphi^u)^2) \boxplus \varphi_0^u$. Denote by
$\theta$ the element of Gal$(\tilde M/\Q)$ which restricts to the
non-trivial automorphism of $K/F$. Appealing to (8.26) and the
compatibility of tensor induction in stages, we get the following
decomposition over $F$:
$$
V_\ell(\pi_f)_F \, \simeq \, {\rm Ind}_L^F(\varphi^2) \otimes {\rm
Ind}_{L^\theta}^F((\varphi^\theta)^2) \oplus \beta_\ell \oplus
\varphi_0\varphi_0^\theta(\nu\nu^\theta \oplus 1 \oplus \nu\delta
\oplus \nu^\theta\delta),
$$
where $\beta_\ell$ is a sum of irreducible $2$-dimensional
representations. If $L \ne L^\theta$, $\nu\delta \oplus
\nu^\theta\delta$ extends to an irreducible of ${\mathcal G}_\Q$,
and moreover, ${\rm Ind}_L^F(\varphi^2) \otimes {\rm
Ind}_{L^\theta}^F((\varphi^\theta)^2)$ is irreducible (already as
a ${\mathcal G}_F$-module). It follows that $r_\ell(\pi_f; \mu)$
is at most $1$ for any Dirichlet character $\mu$, and in this case
we can account for all the Tate classes by algebraic cycles coming
from suitable twists of embedded Hilbert modular surfaces. So
assume that $L = L^\theta$. Then $\nu\nu^\theta = \nu^2 = 1$, and
$r_\ell(\pi_f; \mu)$ is at least $2$ for a suitable $\mu$. {\it
How is one to account for these Tate classes when $K/\Q$ is
cyclic?} Even in the biquadratic case, we can get two independent
algebraic classes only if (a twist of) $\pi$ is a base change all
the way from $\Q$, and things become difficult as seen below.

Now suppose $\pi$ is a base change from $\Q$, i.e., when $\pi =
\pi'_K$ with $\pi' = I_E^\Q(\psi^u)$ for a weight one Hecke
character $\psi$ of an imaginary quadratic field $E$. Let
$\epsilon$, resp. $\nu$, denote the quadratic character of $\Q$
attached to $F$, resp. $E$, and let $\psi'$ denote the transfer of
$\psi$ to ${\mathcal G}_\Q$, which is a finite order character
times $\chi_\ell$. As above let $\delta$ be the quadratic
character of $F$ attached to $K$, with $K/\Q$ Galois. Then we have
${\rm Ind}_F^\Q(\delta) = \gamma \oplus \gamma\epsilon$, where
$\gamma$ is quartic, resp. quadratic, when $K$ is cyclic, resp.
biquadratic. Then, using Proposition 8.28, we see that
$V_\ell(\pi_f)$ is isomorphic to the following:
$$
{\rm Ind}_E^\Q(\psi^4) \oplus \left({\rm Ind}_E^\Q(\psi^2) \otimes
\left(\psi' \oplus \gamma\psi'\epsilon \oplus \gamma\psi' \oplus
\psi'\epsilon\right)\right) \oplus {\psi'}^2\left( 2.{\underline
1} \oplus \nu\gamma\epsilon \oplus \nu\epsilon \oplus \nu\gamma
\oplus \delta_1\right),
$$
where $\delta_1$ is the transfer of $\delta$ to $\Q$. If we write
${\psi'}^2 = \mu\chi_\ell^2$, we then see that $r_\ell(\pi_f;
\mu)$ is $2$ when $K$ is cyclic, and $3$ when $K$ is biquadratic.
So in either case we get an {\it exotic Tate class over an abelian
extension of $\Q$}! Note that when $K$ is biquadratic, since $M$
is the compositum of $K$ with $E$ it is a {\it triquadratic
field}, i.e, $M/\Q$ is Galois with group $(\Z/2)^3$. A natural
question here (because of Lemma 8.3 (ii)) is this: {\it Can such a
$\pi = (\pi')_K$ {\rm (of CM type)} be everywhere unramified with
trivial character?}

Finally suppose we are in the triquadratic case, with $\Psi$ a
weight one Hecke character of $M$ and $\pi = {I}_M^K(\Psi^u)$ {\it
not necessarily a base change from anywhere}. We want to point to
an interesting Tate class. (There are three such, up to complex
conjugation.) Let Gal$(M/\Q)$ be generated by (quadratic elements)
$\rho, \theta, \tau$, with $K$ be the fixed field of (complex
conjugation) $\rho$ and $\theta$ restricting to the non-trivial
automorphism of $K/F$. Put
$$
\Phi : = \, \{1, \rho\theta, \rho\tau, \theta\tau\},
$$
which is a CM type and a subgroup of Gal$(M/\Q)$. If we put
$$
\xi \, = \, \Psi\Psi^{\rho\theta}\Psi^{\rho\tau}\Psi^{\theta\tau},
$$
then it is of Tate type and extends to $\chi_\ell^2$ times a
finite order character $\lambda$, say, of ${\mathcal G}_E$, where
$E$ is the imaginary quadratic field fixed by $\Phi$. Note that
the Tate class is defined over the abelian extension of $E$ cut
out by $\lambda$, and it is defined over an abelian extension of
$\Q$ iff $\lambda$ extends to a Dirichlet character of $\Q$ (which
happens iff $\xi = \xi^\rho$). By starting with a suitable weight
$1$ Hecke character of an imaginary quadratic field and pulling
back to $M$ by norm, it appears that both cases can occur.
D.~Rohrlich has pointed out to the author that one can construct
an {\it everywhere unramified} weight one Hecke character $\Psi$
even when $M/K$ is unramified. {\it Can the corresponding
$\lambda$ be extendable to $\Q$?}

\medskip

\vskip 0.2in

\section{Algebraicity of some Hodge classes}

\bigskip

Now we will prove Theorem B. In view of Propositions 3.6 and 4.11,
it suffices to prove the algebraicity of the Hodge classes in
$V_B(\pi_f)$ for $\pi$ cuspidal. By hypothesis, the level
$\mathfrak N \ne {\mathfrak O}_K$ is square-free. By part (ii) of
Lemma 8.3, if $\pi$ is CM, it will already contrbute at full
level. Put $C = C_0({\mathfrak N})$, $C(1) = C_0({\mathfrak
O}_K)$, $X = \tilde S_C$ and $X(1) = \tilde S_{C(1)}$, with the
associated ramified cover $X \to X(1)$. Since we are interested
only in the Hodge classes on $X_\C$ which are not pull-backs of
classes in $X(1)_\C$, every relevant $\pi$ will be non-CM and its
conductor will be divisible exactly once by a prime $P$. Then, by
part (i) of Lemma 8.3, the local component $\pi_P$ must be an
unramified twist of the Steinberg representation. Fix such a pair
$(\pi,P)$. Choose a quaternion division algebra $B$ over $K$ which
is ramified at $P$ and at three of the infinite places, but
nowhere else. Then by the Jacquet-Langlands correspondence, $\pi$
transfers to a cusp form $\pi'$ on $B^\ast/K$, and $\pi'$
evidently contributes to the $H^1$ of the Shimura curve $C_B$
associated to $B$ (at the corresponding level). This curve is
proper and smooth over $K$, and one knows how to associate (by the
Eichler-Shimura theory, for example) an abelian variety quotient
$A_{\pi}$ of the Jacobian of $C_B$ with the same $L$-function as
that of $\pi$. This leads to the following identity (for a
suitable finite set $S$ of places containing the archimedean
ones):
$$
L^S(s, V_\ell(\pi_f)) \, \simeq \, L^S(s, U_\ell(\pi_f)),
\leqno(9.1)
$$
where
$$
U_\ell(\pi_f) \, = \, \otimes_{\tau} H_\ell^1(A_\pi^\tau
\times_{K^\tau} \overline \Q, \Q_\ell),\leqno(10.2)
$$
where $\tau$ runs over Hom$(K, \R)$ and $A_\pi^\tau$ denotes the
$\tau$-conjugate abelian variety. By Tchebotarev this leads to the
isomorphism as Gal$(\overline \Q/\Q)$-modules:
$$
V_\ell(\pi_f)^{\it ss} \, \simeq \, U_\ell(\pi_f), \leqno(10.3)
$$
where the superscript {\it ss} signifies semi-simplification.

We will now accept the truth of the following result, a special
case of a more general theorem proved in a joint work with
V.~Kumar Murty, which will be published elsewhere:

\medskip

\noindent{\bf Theorem 10.4} \, \it We have
$$
V_B(\pi_f) \, \simeq \, \otimes_{\tau} H_B^1(A_\pi^\tau(\C), \Q),
$$
the isomorphism being one of rational Hodge structures. \rm

\medskip

Note that the right hand sides of the isomorphisms given by (10.3)
and Theorem 10.4 are subspaces of the cohomology (in degree $4$)
of the product abelian variety $\prod_\tau A_\pi^\tau$.

\medskip

By the main theorem of [DMOS] (see section 6 of chapter I), the
$\Q_\ell$-points of the Mumford-Tate group $MT(H^1(A))$ of any
abelian variety $A$ contains the Zariski closure $G(H^1_\ell(A))$
of the image of Galois in the $\ell$-adic representation. By
(10.1), (10.2), the same holds in our case, i.e.,
$MT(V_B(\pi_f))(\Q_\ell) \supset G(V_\ell(\pi_f))$. The Hodge
cycles are none other than the cohomology classes which are fixed
by the Mumford-Tate group; every such class is also fixed by
Galois and hence gives rise to a Tate class. Hence the Tate
conjecture for the $\pi_f$-component of $\tilde S_C$ implies the
Hodge conjecture for this piece. So Theorem B follows from Theorem
A, and we are done.

\vskip 0.2in

\section{Refinements and errata for the papers [Ra2,3]}

\bigskip

\bigskip

In [Ra2] the automorphic tensor product from GL$(2) \times $GL$(2)
\rightarrow $GL$(4)$ was established. In particular one obtained,
given cuspidal automorphic representations $\pi, \pi'$ of GL$(2,
\A_F)$, $F$ any number field, an isobaric automorphic
representation $\pi \boxtimes \pi'$ of GL$(4, \A_F)$ such that
$$
L(s, \pi \boxtimes \pi') \, = \, L(s, \pi \times \pi')
\leqno(11.1)
$$
and
$$
\varepsilon(s, \pi \boxtimes \pi') \, = \, \varepsilon(s, \pi
\times \pi'),
$$
where the functions on the right are respectively the
Rankin-Selberg $L$-function and the associated
$\varepsilon$-factor occurring in its functional equation.

In [Ra2], Theorem M, one also finds a {\it cuspidality criterion},
which is not sufficiently sharp when $\pi$ and $\pi'$ are both
{\it dihedral}. Here is the best possible statement, which should
be used to supplement Theorem M of [Ra2]:

\medskip

\noindent{\bf Theorem 11.2} \, \it Let $\pi, \pi'$ be cuspidal
automorphic representations of GL$(2, \A_F)$. When exactly one of
them is dihedral, $\pi \boxtimes \pi'$ is cuspidal, and when they
are both non-dihedral, $\pi \boxtimes \pi'$ is not cuspidal iff
$\pi'$ is of the form $\pi \otimes \chi$ for an idele class
character $\chi$ of $F$.

Suppose $\pi, \pi'$ are both dihedral, so that we have
$$
\pi = I_E^F(\nu) \quad {\rm and} \quad \pi' = I_{K}^F(\mu),
$$
where $E, K$ are quadratic extensions of $F$, and $\nu$ (resp.
$\nu$) is an idele class character of $E$ (resp. $K$). Then the
following are equivalent:
\begin{enumerate}
\item[($\alpha$)] \, $\pi \boxtimes \pi'$ is not cuspidal
\item[($\beta$)] \, $E = K$
\end{enumerate}
And when one of these equivalent conditions holds, we have
$$
\pi \boxtimes \pi' \, \simeq \, I_K^F(\nu\mu) \boxplus
I_K^F(\nu\mu^\tau),
$$
where $\tau$ is the non-trivial automorphism of $K/F$. \rm

\medskip

{\it Proof of Theorem 11.2} \, When either $\pi$ or $\pi'$ is
non-dihedral, Theorem M of [Ra2] gives what is asserted here. So
assume that both of the representations are dihedral; write them
in the form above. By Theorem M of [Ra2], $\pi \boxtimes \pi'$ is
non-dihedral iff the base change $\pi_K$ is cuspidal and admits a
self-twist by $\lambda: = (\mu \circ \tau)\mu^{-1}$, which is
non-trivial on account of $\pi'$ being cuspidal; we see that
$\lambda$ must be quadratic by comparing the central characters.
Hence $\pi_K$ is cuspidal and dihedral when $\pi \boxtimes \pi'$
is not cuspidal. It suffices to show that $\pi$ is itself dihedral
in such a case. If $L$ is the quadratic extension of $K$ cut out
by $\lambda$, $\pi_K = I_L^K(\nu_L)$ for a character $\nu'$ of
$L$. Then $\pi_K$ corresponds to the (irreducible) 2-dimensional
representation $\rho^K = {\rm Ind}_L^K(\nu')$ of the Weil group
$W_K$. Since $\rho^K$ is invariant under Aut$(K/F)$, it extends to
an irreducible 2-dimensional representation $\rho$ of $W_F$. Since
$\rho$ is of solvable type (and 2-dimensional), it corresponds to
a cusp form $\pi_1$ of GL$(2, \A_F)$ by Langlands and Tunnell.
Moreover, since $\pi_1$ and $\pi$ have the same base change to
$K$, they must differ by at most the quadratic character $\delta$
of $F$ corresponding to K/F. Thus, after replacing $\rho$ by $\rho
\otimes \delta$ if necessary, we may assume that $\rho$ and $\pi$
correspond. It then suffices to show that $\rho$ is dihedral. Let
$\tau$ be the non-trivial automorphism of $K/F$, extended to an
automorphism, again denoted $\tau$, of the Galois closure $\tilde
L$ of $L$ over $F$. Since $\rho^K$ is associated to $L/K$, and
since we have
$$
{\rm End}(\rho_K) \, = \, \rho^K \otimes (\rho^K)^\vee \, = \,
{\rm Ad}(\rho^K) \oplus 1,
$$
its adjoint square Ad$(\rho)$ contains the quadratic character
$\lambda$. There are two cases to consider:

{\bf (i) \, $\lambda = \lambda^\tau$}: \, In this case $\lambda$
extends to a character of $W_F$, and so Ad$(\rho)$ splits as a sum
of such an extension plus a 2-dimensional representation. This
forces $\tau$ to be dihedral. To elaborate, since $\tau$ is
irreducible, Ad$(\rho)$ cannot contain $1$ by Schur's lemma. Since
it is self-dual, if it contains a 1-dimensional summand, it must
contain a quadratic character, say $\epsilon$. Then, if $M/F$ is
the quadratic extension cut out by $\epsilon$, and if $\rho_M$
denotes the restriction $\rho$ to $W_M$, then Ad$(\rho_M)$
contains $1$, hence $\rho_M$ is reducible and $\rho$ must be
induced from $M$.

{\bf (ii) \, $\lambda \ne \lambda^\tau$}: \, Since Ad$(\rho^K)$
extends to $W_F$, it is invariant under $\tau$, and so we are
forced to have the decomposition
$$
{\rm Ad}(\rho^K) \, = \, \lambda \oplus \lambda^\rho \oplus
\lambda\lambda^\tau.
$$
Now the character $\lambda\lambda^\tau$ extends to $W_F$, and
Ad$(\rho)$ again contains a one-dimensional summand, implying that
$\rho$ is dihedral. Done.

\medskip

\noindent{\bf Remark}: \, Erez Lapid has remarked that this can
also be proved by appealing to the properties of the symmetric
square lifting $\pi \to {\rm sym}^2(\pi)$ of Gelbart and Jacquet
from GL$(2)$ to GL$(3)$. More precisely, one appeals to the fact
that $\pi$ is dihedral if and only if sym$^2(\pi)$ is not
cuspidal. Moreover, it can be checked that the Gelbart-Jacquet
lift is compatible with base change, and that quadratic base
change preserves cuspidality for GL$(n)$ for any {\it odd} $n$.
Lapid has had occasion to use Theorem 11.2 in his elegant article
[Lap].

\bigskip

We now move to section 3 of [Ra2], where a key lemma is proved as
a preliminary step to achieving {\bf boundedness in vertical
strips} for the triple product $L$-functions. In fact the key
lemma there is very general and could be of use in various other
situations, giving a bound for the sup norm of an eigenfunction of
an elliptic operator (such as the Laplacian) in terms of its
$L^2$-norm and the eigenvalue $\delta$; in fact it also applies to
functions which are not eigenfunctions. Joseph Shalika recently
asked for a clarification of one of the points of the proof given
in [Ra2], and this what we will do right now. First let us restate
Lemma 3.4.9 as

\medskip

\noindent{\bf Lemma 11.3} \quad \it Let $\Omega$ be a subset of
$\R^N$ contained in the unit ball of radius $\epsilon$, $\Omega'$
a subset of $\Omega$ with non-empty interior such that $\bar
{\Omega'} \, \subset \, \Omega$, and $\Delta$ an elliptic operator
of order $2$. Then we have the following:
\begin{itemize}
\item{(1) For any integer $r > N/2$, there exists a constant
$C_1 > 0$ depending only on
$\epsilon$ and $\Delta$ such that, for all $u$ in the Sobolev
space $\mathcal H_r(\Omega)$,
$$
||u||_{\infty, \Omega} \, \leq \, C_1 \, ||u||_{(2,r); \Omega}.
$$}
\item{(2) For any integer $i \geq 2$, there is a constant $C_i > 0$
depending only on
$\epsilon$ and $\Delta$, such that for any $u \in \mathcal
H_i(\Omega)$:
$$
||u||_{(2,i); \Omega'} \, \leq \, C_i \, \left( ||u||_{2, \Omega}
+ ||\Delta u||_{(2,i-2); \Omega}\right).
$$}
\end{itemize}
\rm

\medskip

Here $||.||_{(2,r); X}$ denotes, for any $r$ and $X = \Omega,
\Omega'$, the $r$th $L^2$-derivative on $X$ so that (for $r \geq 0
$) $||u||_{(2,r); X}$ equals $\sum_{|\nu| \leq r} ||\partial^\nu
u||_{2,X}$, with the $\partial^\nu$ denoting distribution
derivatives. Clearly, $||u||_{(2,0); X} \, = \, ||u||_{2; X}$. We
will henceforth suppress the space $X$ in the subscript.
Furthermore, in Lemma 3.4.9 of [Ra2], $\Delta$ is (at least)
implicitly taken to be the Laplacian, which is not necessary.

The proof of part (1) is as in [Ra2], page 71. Part (2) is proved
by induction on $i \geq 2$. Shalika wanted to know the proof for
the {\bf starting point} $i=2$, and here it is. Choose a nice
cut-off function $\psi$ with compact support in $\Omega$ such that
$\psi u = v$ on $\Omega'$, with $v$ in $\mathcal H_2^0(\Omega)$,
the closure of $C_c^\infty(\Omega)$ in $\mathcal H_r(\Omega)$. By
Theorem 6.2.8 (chapter 6, page 210 of Folland's book [Fo] (second
edition --1995), we have (for all $s \in \R$ and $v \in \mathcal
H_s^0(\Omega')$),
$$
||v||_{(2,s)} \leq C\left(||Lu||_{(2,s-k)} +
||v||_{(2,s-1)}\right),
$$
where $L$ is an elliptic operator of degree $k$ defined on the
closure of $\Omega'$ and $C$ a constant. Applying this with $s=2,
L = \Delta, k=2$ and $v = \psi u$, we get
$$
||\psi u||_{(2,2)} \, \leq \, C\left(||\Delta(\psi u)||_{2} +
||\psi u||_{(2,-1)}\right). \leqno(\ast)
$$
We can evidently bound the left hand side of $(\ast)$ from below
by a constant times $||u||_{(2,1)}$, and also bound $||\psi
u||_{2,-1}$ from above by a constant times $||u||_2$. So to
establish part (2) for $n=2$, we need to bound $||\Delta(\psi
u)||_2$ from above. By the Leibnitz rule we can write $\Delta(\psi
u)$ as $\Delta(\psi)u + \psi\Delta(u)$ plus a sum of terms of the
form $L_1(\Psi)L_2(u)$, where $L_1, L_2$ are differential
operators of order $1$. (When $\Delta$ is the Laplacian, $L_i$ is
the gradient $\nabla$.) Since we can control $L_1(\Psi)$ well, we
need only to bound $||L_2(u)||_{2}$ for any first order operator
$L_2$, or equivalently, we need to bound $||u||_{(2,-1)}$. But
this can be bounded by $||u||_2$ since by the definition of these
spaces using Fourier transform, $||.||_s \leq ||.||_t$ for all
$(s,t)$ with $s \leq t$ (cf. [Fo], page 192). Done.

\medskip

Here are some {\bf typos in this section (3.4) of [Ra2]} which
should be fixed as follows, where $A \, \rightarrow \, B$ means
{\it change $A$ to $B$} and $p.x, \ell.y$ means page $x$, line
$y$, with the understanding that negative line numbers are to be
counted from the bottom of the page:
\begin{enumerate}
\item[] $p.70, \ell.-13$ (Lemma 3.4.8): \, {\it and $\lambda_s$
such that} $\rightarrow$ {\it such that for all $f$} \item[]
$p.71, \ell.9$ (part (2) of Lemma 3.4.9): \, $||u||_{2, \Omega'}
\rightarrow ||u||_{2,\Omega}$ \item[] $p.71, \ell.-9$ : \,
$||\Delta u||_{(2,i-2)} \rightarrow ||\Delta u||_{(2,i-2)}$
\item[] $p.71, \ell.-10$ : \, $||\Delta u||_{(2,i-1)} \rightarrow
||\Delta u||_{(2,i-1)}$ \item[] $p.72, \ell.7$ : \, $\Delta \,
\rightarrow $ the Casimir operator
\end{enumerate}

Further, in the ensuing bound (on page 72 of [Ra2]) of the Arthur
truncation of $E(f_s)$, it should be noted that the Eisenstein
series has a fixed $K$-type.

\bigskip

Now let $K/F$ be a quadratic extension of number fields. In [Ra3]
the {\bf Asai transfer} $\pi \to As(\pi)$ of isobaric automorphic
forms from GL$(2)/K$ to GL$(4)/F$ was achieved. The {\bf sharp
cuspidality criterion} above (Theorem 11.2) has the following {\bf
Asai analogue} in the dihedral case, and this could profitably be
used to replace part (c) of Theorem 1.4 of [Ra3]:

\medskip

\noindent{\bf Theorem 11.4} \, \it Let $K/F$ be a quadratic
extension of number fields with non-trivial automorphism $\theta$.
Let $\pi$ be a dihedral cusp form on GL$(2)/K$, i.e., associated
to a representation $\sigma = I_M^K(\chi)$ of $W_K$ for a
character $\chi$ of a quadratic extension $M$ of $K$. Then the
following are equivalent:
\begin{enumerate}
\item[($\alpha$)] \, $As(\pi)$ is cuspidal
\item[($\beta$)] \, $M$ is not Galois over $F$
\end{enumerate}
And when one of these equivalent conditions fails to hold, there
exist isobaric automorphic representations $\eta, \eta'$ of GL$(2,
\A_F)$ such that
$$
As(\pi) \, \simeq \, \eta \boxplus \eta',
$$
and
$$
\eta_K \, \simeq \, I_M^K(\chi\chi^{\tilde \theta}) \quad {\rm
and} \quad \eta'_K \, \simeq \, I_M^K(\chi\chi^{\tilde
\theta\tau}),
$$
where $\tau$ is the non-trivial automorphism of $M/K$ and $\tilde
\theta$ denotes an extension of $\theta$ to an automorphism of $M$
over $F$. \rm

\medskip

{\bf Remark 11.5}: \, When $\pi$ is non-dihedral (and cuspidal),
as proved in part (b) of Theorem 1.4 of [Ra3], $As(\pi)$ is a
cuspidal automorphic representation of GL$(4, \A_F)$ iff $\pi
\circ \theta$ is not isomorphic to $\pi \otimes \mu$ for any idele
class character $\mu$ of $K$. There is no refinement in this case.

\medskip

{\it Proof}. \, Recall from [Ra3] that the base change $As(\pi)_K$
of $As(\pi)$ to GL$(4)/K$ satisfies
$$
As(\pi)_K \, \simeq \, \pi \boxtimes (\pi \circ \theta).
\leqno(11.6)
$$
Since $\pi$ corresponds to $\sigma = {\rm Ind}_M^K(\chi)$, this
translates to the following isomorphism of $W_K$-modules:
$$
As(\sigma)_K \, \simeq \, \sigma \boxtimes \sigma^{[\theta]}.
\leqno(11.7)
$$

When $M$ is non-Galois, Theorem 11.2 implies that $\pi \boxtimes
(\pi \circ \theta)$ is cuspidal. By (11.6), this is the base
change of $As(\pi)$ to GL$(4)/K$. It follows that $As(\pi)$ is
itself cuspidal. Hence ($\beta$) implies ($\alpha$).

It was proved in [Ra3] that ($\alpha$) implies ($\beta$). Here is
a slightly different proof. Suppose $M$ is Galois over $F$, i.e.,
$M = M^{\tilde \theta}$. Then we get the following decomposition
by Mackey and (11.7):
$$
As(\sigma)_K \, \simeq \, {\rm Ind}_M^K(\chi\chi^{\tilde \theta})
\oplus {\rm Ind}_M^K(\chi\chi^{{\tilde \theta}\tau}) \leqno(11.8)
$$
The first representation on the right is evidently invariant under
$\theta$ and consequently extends to a $2$-dimensional
representation $\eta$, say, of $W_F$. This means that $As(\sigma)$
is reducible, and so the function  $L(s, As(\sigma) \otimes
As(\sigma)^\vee)$ has at least a double pole at $s=1$. Since by
[Ra3] this $L$-function is the same as the Rankin-Selberg
$L$-function $L(s, As(\pi) \times As(\pi)^\vee)$, this automorphic
$L$-function also has a pole of order $\geq 2$ at $s=1$. But then
the results of Jacquet and Shalika imply that $As(\pi)$ is {\it
not} cuspidal. Done.

\bigskip

Next we need to make a {\bf correction of incompatible sign
conventions used in [Ra3]}. This is important, even though the
incompatibilities did not affect the truth of any of the main
results in that paper. To elaborate, let $K/F$ be a quadratic
extension with non-trivial automorphism $\theta$, and let $\delta$
be the quadratic character of ${\mathcal G_F}$ associated to
$K/F$. Given an irreducible $2$-dimensional representation
$\sigma$ of ${\mathcal G}_K$) (or $W_K$), the Asai representation
of $\sigma$ is a choice of an extension $A(\sigma)$, say, to
${\mathcal G}_F$ of the $\theta$-invariant representation $\sigma
\otimes \sigma^{[\theta]}$. There are at least two such extensions
as one can replace $A(\sigma)$ by its twist by the {\it sign
character} $\delta$; these are the only extensions when
$A(\sigma)$ is irreducible. On page 19 (see (4.50)) of [Ra3], we
defined the Asai representation to be the summand of
$\Lambda^2({\rm Ind}_K^F(\sigma))$ with complement ${\rm
Ind}_K^F({\rm det}(\sigma))$. This, as well as the analogous
definition of $As(\pi)$ for a cuspidal $\pi$ on GL$(2)/K$, is fine
till we get to the asserted identity (4.64) on page 21 (of [Ra3]),
which is off by the sign character $\delta$; we need to twist
(exactly) one of the $L$-functions by $\delta$. The reason is
this: Langlands's definition of $r$ is compatible with {\it tensor
induction}, which really occurs in the {\it symmetric square} of
the usual induction; the $\delta$-twist of it occurs in the
exterior square. It perhaps makes sense to define $As(\sigma)$
(and $As(\pi)$) using tensor induction. Then we should implement
the following {\bf errata} to [Ra3]:

\begin{enumerate}
\item[] $p.19, \ell.-4, \, (4.50)$: \, $As(\sigma) \, \rightarrow \, As(\sigma)
\otimes \delta$
\item[] $p.31, \ell.-11, \, (6.4)$: \,
$$
{\rm sym}^2(\pi_0) \otimes \delta(\mu\nu)^{-1} \boxplus
\delta\mu^{-1}\, \rightarrow \, {\rm sym}^2(\pi_0) \otimes
(\mu\nu)^{-1} \boxplus \mu^{-1}
$$
\item[] $p.32, \ell.7$: \, the induction \, $\rightarrow$ \, the
$(\nu_v\delta_v)^{-1}$-twist of the induction
\item[]$p.32, \ell.9, \, (6.8)$: \,
$$
{\rm sym}^2(\sigma_v(\pi_0)) \otimes \delta_v(\mu_v\nu_v)^{-1}
\oplus \delta_v\mu_v^{-1} \, \rightarrow \, {\rm
sym}^2(\sigma_v(\pi_0)) \otimes (\mu_v\nu_v)^{-1} \oplus
\mu_v^{-1}
$$
\item[]$p.37, \ell.-4, \, (6.8)$: \, $\mu_{1,0}\delta_v \oplus
\mu_{2,0}\delta_v \, \rightarrow \, \mu_{1,0} \oplus \mu_{2,0}$
\item[]$p.37, \ell.-2, \, (6.8)$: \, $\mu_{1,0}\delta_v \boxplus
\mu_{2,0}\delta_v \, \rightarrow \, \mu_{1,0} \boxplus \mu_{2,0}$
\end{enumerate}

\medskip

Finally, it is remarked on page 22 of [Ra3] that the Asai transfer
$\Pi$ to GL$(4)/F$ satisfies $L(s, \Pi) = L(s, \pi; r)$. But in
fact the method of the paper implies more (see sections 7, 8),
namely that for any isobaric automorphic representation $\eta$ of
GL$(2, \A_F)$,
$$
L(s, \Pi \times \eta) \, = \, L(s, \pi; r \otimes \eta).
$$
In particular, the local component $\Pi_v$ is at any place $v$ the
correct admissible representation of GL$(4, F_v)$ associated by
functoriality.

\vskip 0.2in

\addcontentsline{toc}{section}{Bibliography}
\section*{\bf Bibliography}

\begin{description}

\item[{[AC]}] J.~Arthur and L.~Clozel, {\it Simple Algebras, Base Change
and the Advanced Theory of the Trace Formula}, Ann. Math. Studies {\bf 120}
(1989), Princeton, NJ.

\item[{[As]}] T.~Asai, On certain Dirichlet series associated with Hilbert modular
forms and Rankin's method,  Math. Annalen  {\bf 226}, no. 1
(1977), 81--94.

\item[{[BBD]}] A.~Beilinson, J.~Bernstein and P.~Deligne,
Faisceaux pervers, in {\it Analysis and topology on singular
spaces I},  5--171, As\'erisque, 100, Soc. Math. France, Paris
(1982).

\item[{[B$\ell$]}] D.~Blasius, Hilbert modular forms and the Ramanujan conjecture,
preprint (2003).

\item[{[B$\ell$-Ro]}] D.~Blasius and J.~Rogawski, Motives for
Hilbert modular forms, Invent. Math. {\bf 114} (1993), no. 1,
55--87.

\item[{[BoJ]}] A.~Borel and H.~Jacquet, Automorphic forms and
automorphic representations, in {\it Automorphic forms,
representations and $L$-functions}, Proc. Symposia Pure Math.,
XXXIII, Part 1, 189--207, AMS, Providence, R.I. (1979).

\item[{[BoW]}] A.~Borel and N.~Wallach, {\it Continuous
cohomology, discrete subgroups, and representations of reductive
groups}, second edition, Mathematical Surveys and Monographs {\bf
67}, AMS, Providence, RI (2000).

\item[{[BrL]}] J.-L.~Labesse, Cohomologie d'intersection et
fonctions $L$ de certaines variét\'es de Shimura, Ann. Sci. École
Norm. Sup. (4) 17 (1984), no. 3, 361--412.

\item[{[CuR]}] C.W.~Curtis and I.~Reiner, Methods of
representation theory I, Wiley, NY (1981).

\item[{[De1]}] P.~Deligne,  Travaux de Shimura, S\'eminaire
Bourbaki, 23ème ann\'ee (1970/71), Expos\'e {\it 389}, 123--165,
Lecture Notes in Math. {\bf 244}, Springer, Berlin (1971).

\item[{[De2]}] P.~Deligne, Les constantes des équations
fonctionnelles des fonctions $L$, in {\it Modular functions of one
variable II}, 501--597, Lecture Notes in Math. {\bf 349},
Springer, Berlin (1973).

\item[{[DMOS]}] P.~Deligne, J.~Milne, A.~Ogus and K.~Shih, {\it
Hodge cycles, motives, and Shimura varieties}, Lecture Notes in
Math. {\bf 900}, Springer-Verlag, NY, 1982.

\item[{[F]}] G.~Faltings,  $p$-adic Hodge theory, Journal of the
AMS {\bf 1} (1988), no. 1, 255--299.

\item[{[F$\ell$]}] Y.Z.~Flicker, Twisted tensors and Euler
products {\bf 116} (1988), 295--313.

\item[{[Fo]}] Folland, G.~Folland, {\it Introduction to partial
differential equations}, Second edition, Princeton University
Press, Princeton, NJ (1995).

\item[{[Ge]}] S.~Gelbert, {\it Automorphic forms on adele groups},
Annals of Math. Studies {\bf 83}, Princeton, NJ (1975).

\item[{[GeJ]}] S.~Gelbart and H.~Jacquet, A relation between
automorphic representations of GL$(2)$ and GL$(3)$, Ann. Scient.
\'Ec. Norm. Sup. (4) {\bf 11} (1979), 471--542.

\item[{[HLR]}] G.~Harder, R.P.~Langlands and
M.~Rapoport, Algebraische zykeln auf Hilbert-Blumenthal-Fl\"achen,
Crelles Journal {\bf 366} (1986), 53--120.

\item[{[HaT]}] M.~Harris and R.~Taylor, On the geometry and cohomology of some
simple Shimura varieties, Annals of Math. Studies, Princeton
(2001).

\item[{[He]}] G.~Henniart, Une preuve simple des conjectures de Langlands pour
${\rm GL}(n)$ sur un corps $p$-adique, Invent. Math. {\bf 139}, no. 2, 439--455
(2000).

\item[{[HZ]}] F.~Hirzebruch and D.~Zagier, Intersection numbers of curves on
Hilbert modular surfaces and modular forms of Nebentypus, Invent.
Math. {\bf 36} (1976), 57--113.

\item[{[Ik]}] T.~Ikeda, On the location of poles of the triple
$L$-functions,  Compositio Math. {\bf 83} (1992), no. 2, 187--237.

\item[{[JPSS]}] H.~Jacquet, I.~Piatetski-Shapiro and J.~Shalika,
Rankin-Selberg convolutions, Amer. J. Math.{\bf 105}, 367--464
(983).

\item[{[JS]}] H.~Jacquet and J.A.~Shalika, On Euler products and
the classification of automorphic forms I \& II, Amer. J of Math.
{\bf 103} (1981), 499--558 \& 777--815.

\item[{[JY]}] H.~Jacquet and Y.~Ye,  Une remarque sur le
changement de base quadratique, C. R. Acad. Sci. Paris S\'er. I
Math. {\bf 311} (1990), no. 11, 671--676..

\item[{[K]}] H.~Kim, Functoriality for the exterior square of GL$_4$ and the
symmetric fourth power of GL$_2$, Journal of the AMS {\bf 16}, No.
1 (2002), 139--183.

\item[{[KSh]}] H.~Kim and F.~Shahidi, Functorial products for GL$_2 \times $GL$_3$
and functorial symmetric cube for GL$_2$, with an appendix by
C.J.~Bushnell and G.~Henniart,  Ann. of Math. (2)  {\bf 155}
(2002), no. 3, 837--893.

\item[{[Kr]}] M.~Krishnamurthy, The Asai transfer to GL4 via the
Langlands-Shahidi method, IMRN {\bf 2003}, no. 41, 2221--2254.

\item[{[La1]}] R.P.~Langlands, On the zeta functions of some
simple Shimura varieties,  Canadian J. Math. {\bf 31} (1979), no.
6, 1121--1216.

\item[{[La2]}] R.P.~Langlands, On the notion of an automorphic representation.
A supplement, in {\it Automorphic forms, Representations and
$L$-functions}, ed. by A. Borel and W. Casselman, Proc. symp. Pure
Math {\bf 33}, part 1, 203--207, AMS. Providence (1979).

\item[{[La3]}] R.P.~Langlands, On the classification of
irreducible representations of real algebraic groups, in {\it
Representation theory and harmonic analysis on semisimple Lie
groups}, 101--170, Math. Surveys Monographs {\bf 31}, AMS,
Providence, RI (1989).

\item[{[Lap]}] E.~Lapid,  On the nonnegativity of Rankin-Selberg
$L$-functions at the center of symmetry, IMRN {\bf 2003}, no. 2,
65--75.

\item[{[Lo]}] E.~Looijenga,  $L\sp 2$-cohomology of locally
symmetric varieties, Compositio Math. {\bf 67} (1988), no. 1,
3--20.

\item[{[MuP]}] V.K.~Murty and D.~Prasad, Tate cycles on a product
of Hilbert modular surfaces, Journal of Number theory {\bf 80}
(200), 25--43.

\item[{[MuRa]}] V.K.~Murty and D.~Ramakrishnan, Period relations
and the Tate conjecture for Hilbert modular surfaces, Inventiones
Math. {\bf 89}, no. 2 (1987), 319--345.

\item[{[N]}] A.~Nair, Intersection cohomology, Shimura varieties,
and motives, preprint (2003).

\item[{[PS-R]}] I. Piatetski-Shapiro and S. Rallis, Rankin triple
$L$-functions, Compositio Mathematica {\bf 64} (1987), 31-115.

\item[{[Ra1]}] D.~Ramakrishnan, Arithmetic of Hilbert-Blumenthal surfaces, CMS Conference
Proceedings {\bf 7}, 285--370 (1987).

\item[{[Ra2]}] D.~Ramakrishnan, Modularity of the Rankin-Selberg $L$-series, and multiplicity
one for SL$(2)$, Annals of Mathematics {\bf 152}, 43--108 (2000).

\item[{[Ra3]}] D.~Ramakrishnan, Modularity of solvable Artin
representations of GO$(4)$-type, International Mathematics
Research Notices (IMRN) {\bf 2002}, No. 1 (2002), 1--54.

\item[{[SaSt]}] L.~Saper and M.~Stern, $L\sb 2$-cohomology of
arithmetic varieties, Ann. of Math. (2) {\bf 132} (1990), no. 1,
1--69.

\item[{[Sh]}] F.~Shahidi, On the Ramanujan conjecture and
finiteness of poles for certain $L$-functions.  Ann. of Math. (2)
{\bf 127} (1988), no. 3, 547--584.

\item[{[Ta]}] R.L.~Taylor, On Galois representations associated to Hilbert modular forms,
Invent. Math. {\bf 98} (1989), no. 2, 265--280.

\bigskip

\end{description}

\vskip 0.3in

\noindent
Dinakar Ramakrishnan

\noindent
Department of Mathematics

\noindent
California Institute of Technology, Pasadena, CA 91125.

\noindent
dinakar@its.caltech.edu

\end{document}